\documentclass{article}

\usepackage{natbib}
\usepackage{amsmath,amsthm,amsopn,amstext,amscd,amsfonts,amssymb}
\usepackage{graphicx}

\def \tr {\mathop{\rm tr}\nolimits}
\def \re {\mathop{\rm Re}\nolimits}
\def \im {\mathop{\rm Im}\nolimits}
\def \Tec {\mathop{\rm vec}\nolimits}
\def \Vol {\mathop{\rm Vol}\nolimits}

\def \cov {\mathop{\rm Cov}\nolimits}
\def \etr {\mathop{\rm etr}\nolimits}
\def \diag {\mathop{\rm diag}\nolimits}
\def \build#1#2#3{\mathrel{\mathop{#1}\limits^{#2}_{#3}}}

\newtheorem{theorem}{Theorem}[section]
\newtheorem{lemma}[theorem]{Lemma}
\newtheorem{corollary}[theorem]{Corollary}

\theoremstyle{definition}
\newtheorem{definition}[theorem]{Definition}

\theoremstyle{remark}
\newtheorem{remark}[theorem]{Remark}

\numberwithin{equation}{section}

\renewenvironment{abstract}
                 {\vspace{6pt}
                  \begin{center}
                  \begin{minipage}{5in}
                  \centerline{\textbf{Abstract}}
                  \noindent\ignorespaces
                 }
                 {\end{minipage}\end{center}}

\title{\Large \textbf{Special Functions:\\ Integral properties of Jack polynomials,
    hypergeometric functions and invariant polynomials}}
\author{
  \textbf{Jos\'e A. D\'{\i}az-Garc\'{\i}a} \thanks{Corresponding author\newline
   {\bf Key words.}  Jack polynomials, invariant polynomials, zonal polynomials
    generalised beta and gamma functions, generalised hypergeometric functions
   , real, complex, quaternion and octonion random  matrices.\newline
    2000 Mathematical Subject Classification. Primary 62H10; 15A52; 43A85; 43A90; 32A07; secondary
    32M15; 44A10; 44A15; 62E15}\\
  {\normalsize Department of Statistics and Computation} \\
  {\normalsize 25350 Buenavista, Saltillo, Coahuila, Mexico} \\
  {\normalsize E-mail: jadiaz@uaaan.mx} \\
}
\date{}
\begin{document}
\maketitle

\begin{abstract}
Some integral properties of Jack polynomials, hypergeometric functions and invariant
polynomials are studied for real normed division algebras.
\end{abstract}

\section{Introduction}\label{sec1}

During the 1960s, real and complex zonal polynomials were studied exhaustively by
\citet{j:61, j:64}, \citet{c:63} and \citet{k:66}, among many others. Excellent
reference books include those by \citet{m:82}, \citet{t:84}, \citet{f:85} and
\citet{m:97}, which summarise many of the results published to date.

Hypergeometric functions with a matrix argument were first studied by \citet{h:55}
and defined in terms of zonal  polynomials by \citet{c:63}. Hypergeometric functions
of one or two matrix arguments have been applied in many areas of science and
technology, including multivariate statistical analysis (\citet{m:82} and
\citet{m:97}), random matrix theory (\citet{me:91} and \citet{f:05}), wireless
communications (\citet{rva:05a, rva:05b}, shape theory \citet{gm:93} and
\citet{cldggf:09}).

Later \citet{da:79}, \citet{da:80}, \citet{ch:80} and \citet{chd:86}, introduced a
class of homogeneous invariant polynomials with two or more matrix arguments, which
generalise the zonal polynomials; many of their basic and integral properties are
studied in real cases.

In the context of multivariate statistics, zonal polynomials were initially used to
express many noncentral matrix variate distributions. However, there were other
distributional problems that could not be solved using zonal polynomials. In these
latter cases, invariant polynomials were used to obtain explicit expressions of
doubly noncentral matrix variate distributions, matrix variate distribution functions
and the joint density of eigenvalues of matrix variate beta type I and II
distributions, etc. see \citet{j:64} and \citet{da:80}.

During the 1990s, zonal polynomials regained prominence but from a more general point
of view, in which it was observed that zonal polynomial for real and complex cases
are particular cases of Jack polynomials, see \citet{S97} and \citet{gj96}, among
many others. In terms of Jack polynomials, it is possible to give a general
definition for hypergeometric functions, see \citet{gr:87}, \citet{gr:89}, among many
others. In both cases, Jack polynomials and generalised hypergeometric functions are
written in term of a parameter, denoted by $\alpha$ or $\beta$ and with which, for
example for $\beta = 1$, $2$ o $4$, the zonal polynomials and hypergeometric
functions  are obtained for real, complex and quaternion cases, respectively, see
\citet{gr:87}, \citet{mopsj05} and \citet{ke:06}.

The properties for Jack polynomials and hypergeometric functions with a matrix
argument have been studied by \citet{h:55}, \citet{j:64}, \citet{c:63}, \citet{k:66}
and \citet{m:82} in the real case (zonal polynomials); by \citet{j:64}, \citet{t:84},
\citet{f:85} and \citet{rva:05a,rva:05b} in the complex case (Schur functions); by
\citet{lx:09} in the quaternion case and by \citet{gr:87} and \citet{cldggf:07} in
the general case (real, complex and quaternion cases), among many others.

A serious obstacle encountered when Jack polynomials, hypergeometric functions and
invariant polynomials are to be used is the question of their calculation.
Fortunately, with the excellent algorithm proposed by \citet{ke:06} and \citet{k:04},
it is now possible to use these techniques in many applications, see \citet{rva:05b}
and \citet{cldggf:09}. Unfortunately, this obstacle remains for the general case of
invariant polynomials.

Let us take into account that there are exactly four normed division algebras: the
real numbers ($\Re$), complex numbers ($\mathfrak{C}$), quaternions ($\mathfrak{H}$),
and octonions ($\mathfrak{O}$); moreover, these are the only alternative division
algebras, and all division algebras have a real dimension of $1, 2, 4$ or $8$, see
\citet[Theorems 1, 2 and 3]{b:02}. Furthermore, according to \citet{b:02}, there is
still no proof that octonions are useful for understanding the real world.

In this paper, we generalise diverse integral properties of Jack polynomials,
hypergeometric functions and invariant polynomials for normed division algebras. Note
that we can only conjecture the results for the octonion case, because many of its
related matrix problems are still under study. However, for example in \citet[Section
1.4.5, pp. 22-24]{f:05} it is proved that the bi-dimensional density function of the
eigenvalue, for a $2 \times 2$ octonionic matrix with symmetric normal distribution,
is obtained from the general joint density function of the eigenvalues for the
symmetric normal distribution, assuming $m = 2$ and $\beta = 8$. The material in the
present paper is organised as follows: Section \ref{sec2} provides some notation and
preliminary results about Jacobians, gamma and beta multivariate functions and
invariant measures. The definition and many integral properties of Jack polynomials
are obtained in Section \ref{sec3}. Many extensions of the integral properties of
hypergeometric functions with one and two arguments are studied in Section
\ref{sec4}. For invariant polynomials with two matrix arguments, in Section
\ref{sec5} we derive diverse integral properties, such as the inverse Laplace
transformation, gamma and beta integrals, etc. Finally, in Section \ref{sec6}, we
show diverse applications of some results derived previously, such as the
distribution function of a central Wishart distribution for normed division algebras,
its joint eigenvalue density and the distribution function of the largest and
smallest eigenvalues. We emphasise the conditions that must be met by the parameters
that take part in many integral properties in the cases discussed, because, even in
the original references, these conditions were omitted or established inexactly.

\section{Preliminary}\label{sec2}

Let ${\mathcal L}^{\beta}_{m,n}$ be the linear space of all $n \times m$ matrices of
rank $m \leq n$ over $\mathfrak{F}$ with $m$ distinct positive singular values, where
$\mathfrak{F}$ denotes a \emph{real finite-dimensional normed division algebra}. In
particular, let $GL(m,\mathfrak{F})$ be the space of all invertible $m \times m$
matrices over $\mathfrak{F}$. Let $\mathfrak{F}^{n \times m}$ be the set of all $n
\times m$ matrices over $\mathfrak{F}$. The dimension of $\mathfrak{F}^{n \times m}$
over $\Re$ is $\beta mn$. And let us recall that the parameter $\beta$ has
traditionally been used to count the real dimension of the underlying normed division
algebra. In other branches of mathematics, the parameter $\alpha = 2/\beta$ is used,
see \citet{b:02}.

\medskip

\begin{table}[!h]
  \caption{Values of the $\beta = 2/\alpha$ parameter.}\label{table1}
  \renewcommand\arraystretch{1.5}
\begin{center}
\begin{tabular}{l|l l}
  \hline\hline
  $\beta$ & $\alpha$ & Normed divison algebra \\
  \hline
  1 & 2 & real ($\Re$) \\
  2 & 1 & complex ($\mathfrak{C}$) \\
  4 & 1/2 & quaternion ($\mathfrak{H}$) \\
  8 & 1/4 & octonion ($\mathfrak{O}$) \\
  \hline
\end{tabular}
\end{center}
\end{table}

Let $\mathbf{A} \in \mathfrak{F}^{n \times m}$, then $\mathbf{A}^{*} =
\overline{\mathbf{A}}^{T}$ denotes the usual conjugate transpose. The set of matrices
$\mathbf{H}_{1} \in \mathfrak{F}^{n \times m}$ such that
$\mathbf{H}_{1}^{*}\mathbf{H}_{1} = \mathbf{I}_{m}$ is a manifold denoted ${\mathcal
V}_{m,n}^{\beta}$, termed the \emph{Stiefel manifold} ($\mathbf{H}_{1}$ are also
known as \emph{semi-orthogonal} ($\beta = 1$), \emph{semi-unitary} ($\beta = 2$),
\emph{semi-symplectic} ($\beta = 4$) and \emph{semi-exceptional type} ($\beta = 8$)
matrices, see \citet{dm:99}). The dimension of $\mathcal{V}_{m,n}^{\beta}$ over $\Re$
is $[\beta mn - m(m-1)\beta/2 -m]$. In particular, ${\mathcal V}_{m,m}^{\beta}$ with
dimension over $\Re$, $[m(m+1)\beta/2 - m]$, is the maximal compact subgroup
$\mathfrak{U}^{\beta}(m)$ of ${\mathcal L}^{\beta}_{m,m}$ and consists of all
matrices $\mathbf{H} \in \mathfrak{F}^{m \times m}$ such that
$\mathbf{H}^{*}\mathbf{H} = \mathbf{I}_{m}$. Therefore, $\mathfrak{U}^{\beta}(m)$ is
the \emph{real orthogonal group} $\mathcal{O}(m)$ ($\beta = 1$), the \emph{unitary
group} $\mathcal{U}(m)$ ($\beta = 2$), the \emph{compact symplectic group}
$\mathcal{S}p(m)$ ($\beta = 4$) or \emph{exceptional type matrices} $\mathcal{O}o(m)$
($\beta = 8$), for $\mathfrak{F} = \Re$, $\mathfrak{C}$, $\mathfrak{H}$ or
$\mathfrak{O}$, respectively. Denote by ${\mathfrak S}_{m}^{\beta}$ the real vector
space of all $\mathbf{S} \in \mathfrak{F}^{m \times m}$ such that $\mathbf{S} =
\mathbf{S}^{*}$. Let $\mathfrak{P}_{m}^{\beta}$ be the \emph{cone of positive
definite matrices} $\mathbf{S} \in \mathfrak{F}^{m \times m}$; then
$\mathfrak{P}_{m}^{\beta}$ is an open subset of ${\mathfrak S}_{m}^{\beta}$. Over
$\Re$, ${\mathfrak S}_{m}^{\beta}$ consist of \emph{symmetric} matrices; over
$\mathfrak{C}$, \emph{Hermitian} matrices; over $\mathfrak{H}$, \emph{quaternionic
Hermitian} matrices (also termed \emph{self-dual matrices}) and over $\mathfrak{O}$,
\emph{octonionic Hermitian} matrices. Generically, the elements of
$\mathfrak{S}_{m}^{\beta}$ are termed as \textbf{Hermitian matrices}, irrespective of
the nature of $\mathfrak{F}$. The dimension of $\mathfrak{S}_{m}^{\beta}$ over $\Re$
is $[m(m-1)\beta+2]/2$. Let $\mathfrak{D}_{m}^{\beta}$ be the \emph{diagonal
subgroup} of $\mathcal{L}_{m,m}^{\beta}$ consisting of all $\mathbf{D} \in
\mathfrak{F}^{m \times m}$, $\mathbf{D} = \diag(d_{1}, \dots,d_{m})$. Let
$\mathfrak{T}_{U}^{\beta}(m)$ be the subgroup of all \emph{upper triangular} matrices
$\mathbf{T} \in \mathfrak{F}^{m \times m}$ such that $t_{ij} = 0$ for $1 < i < j \leq
m$; and let $\mathfrak{T}_{L}^{\beta}(m)$ be the opposed \emph{lower triangular}
subgroup $\mathfrak{T}_{L}^{\beta}(m) =
\left(\mathfrak{T}_{U}^{\beta}(m)\right)^{T}$. For any matrix $\mathbf{X} \in
\mathfrak{F}^{n \times m}$, $d\mathbf{X}$ denotes the\emph{ matrix of differentials}
$(dx_{ij})$. Finally, we define the \emph{measure} or volume element $(d\mathbf{X})$
when $\mathbf{X} \in \mathfrak{F}^{m \times n}, \mathfrak{S}_{m}^{\beta}$,
$\mathfrak{D}_{m}^{\beta}$ or $\mathcal{V}_{m,n}^{\beta}$, see \citet{d:02}.

If $\mathbf{X} \in \mathfrak{F}^{n \times m}$ then $(d\mathbf{X})$ (the Lebesgue
measure in $\mathfrak{F}^{n \times m}$) denotes the exterior product of the $\beta
mn$ functionally independent variables
$$
  (d\mathbf{X}) = \bigwedge_{i = 1}^{n}\bigwedge_{j = 1}^{m}\bigwedge_{k =
  1}^{\beta}dx_{ij}^{(k)}.
$$
\begin{remark}
\begin{itemize} Note that for $x_{ij} \in \mathfrak{F}$
$$
    dx_{ij} = \bigwedge_{k = 1}^{\beta}dx_{ij}^{(k)}.
  $$
In particular for $\mathfrak{F} = \Re$, $\mathfrak{C}$, $\mathfrak{H}$ or
$\mathfrak{O}$ we have
  \item $x_{ij} \in \Re$ then
  $$
    dx_{ij} = \bigwedge_{k = 1}^{1}dx_{ij}^{(k)}= dx_{ij}.
  $$
  \item $x_{ij} = x_{ij}^{(1)} + ix_{ij}^{(2)} \in \mathfrak{C}$, then
  $$
    dx_{ij} = dx_{ij}^{(1)} \wedge dx_{ij}^{(2)} = \bigwedge_{k = 1}^{2}dx_{ij}^{(k)}.
  $$
  \item $x_{ij} = x_{ij}^{(1)} + ix_{ij}^{(2)}+ jx_{ij}^{(3)}+ kx_{ij}^{(4)} \in
  \mathfrak{H}$, then
  $$
    dx_{ij} = dx_{ij}^{(1)} \wedge dx_{ij}^{(2)}\wedge dx_{ij}^{(3)}\wedge dx_{ij}^{(4)}
    = \bigwedge_{k = 1}^{4}dx_{ij}^{(k)}.
  $$
  \item $x_{ij} = x_{ij}^{(1)} + e_{1}x_{ij}^{(2)}+ e_{2}x_{ij}^{(3)}+ e_{3}x_{ij}^{(4)} +
  e_{4}x_{ij}^{(5)} + e_{5}x_{ij}^{(6)}+ e_{6}x_{ij}^{(7)}+ e_{7}x_{ij}^{(8)} \in
  \mathfrak{O}$, then
  $$
    dx_{ij} = dx_{ij}^{(1)} \wedge dx_{ij}^{(2)} \wedge dx_{ij}^{(3)} \wedge dx_{ij}^{(4)}
    \wedge dx_{ij}^{(5)} \wedge dx_{ij}^{(6)} \wedge dx_{ij}^{(7)} \wedge
    dx_{ij}^{(8)} = \bigwedge_{k = 1}^{8}dx_{ij}^{(k)}
  $$
\end{itemize}
\end{remark}
If $\mathbf{S} \in \mathfrak{S}_{m}^{\beta}$ (or $\mathbf{S} \in
\mathfrak{T}_{L}^{\beta}(m)$) then $(d\mathbf{S})$ (the Lebesgue measure in
$\mathfrak{S}_{m}^{\beta}$ or in $\mathfrak{T}_{L}^{\beta}(m)$) denotes the exterior
product of the $m(m+1)\beta/2$ functionally independent variables (or denotes the
exterior product of the $m(m-1)\beta/2 + n$ functionally independent variables, if
$s_{ii} \in \Re$ for all $i = 1, \dots, m$)
$$
  (d\mathbf{S}) = \left\{
                    \begin{array}{ll}
                      \displaystyle\bigwedge_{i \leq j}^{m}\bigwedge_{k = 1}^{\beta}ds_{ij}^{(k)}, &  \\
                      \displaystyle\bigwedge_{i=1}^{m} ds_{ii}\bigwedge_{i < j}^{m}\bigwedge_{k = 1}^{\beta}ds_{ij}^{(k)}, &
                       \hbox{if } s_{ii} \in \Re.
                    \end{array}
                  \right.
$$
\begin{remark}
Since generally the context establishes the conditions on the elements of
$\mathbf{S}$, that is, if $s_{ij} \in \Re$, $\in \mathfrak{C}$, $\in \mathfrak{H}$ or
$ \in \mathfrak{O}$. It shall be considered
$$
  (d\mathbf{S}) = \bigwedge_{i \leq j}^{m}\bigwedge_{k = 1}^{\beta}ds_{ij}^{(k)}
   \equiv \bigwedge_{i=1}^{m} ds_{ii}\bigwedge_{i < j}^{m}\bigwedge_{k =
1}^{\beta}ds_{ij}^{(k)}.
$$
Observe, too, that for the Lebesgue measure $(d\mathbf{S})$ defined thus, it is
required that $\mathbf{S} \in \mathfrak{P}_{m}^{\beta}$, that is, $\mathbf{S}$ must
be a non singular Hermitian matrix (Hermitian definite positive matrix). In the real
case, when $\mathbf{S}$ is a positive semidefinite matrix, its corresponding measure
is studied in  \citet{u:94}, \citet{dg:97}, \citet{dggf:04a} and \citet{dggf:04b}
under different coordinate systems.
\end{remark}
If $\mathbf{\Lambda} \in \mathfrak{D}_{m}^{\beta}$ then $(d\mathbf{\Lambda})$ (the
Legesgue measure in $\mathfrak{D}_{m}^{\beta}$) denotes the exterior product of the
$\beta m$ functionally independent variables
$$
  (d\mathbf{\Lambda}) = \bigwedge_{i = 1}^{n}\bigwedge_{k = 1}^{\beta}d\lambda_{i}^{(k)}.
$$
If $\mathbf{H}_{1} \in \mathcal{V}_{m,n}^{\beta}$ then
$$
  (\mathbf{H}^{*}_{1}d\mathbf{H}_{1}) = \bigwedge_{i=1}^{n} \bigwedge_{j =i+1}^{m}
  \mathbf{h}_{j}^{*}d\mathbf{h}_{i}.
$$
where $\mathbf{H} = (\mathbf{H}_{1}|\mathbf{H}_{2}) = (\mathbf{h}_{1}, \dots,
\mathbf{h}_{m}|\mathbf{h}_{m+1}, \dots, \mathbf{h}_{n}) \in \mathfrak{U}^{\beta}(m)$.
It can be proved that this differential form does not depend on the choice of the
matrix $\mathbf{H}_{2}$ and that it is invariant under the transformations
\begin{equation}\label{unif}
    \mathbf{H}_{1} \rightarrow \mathbf{QHP}_{1},  \quad \mathbf{Q }\in
    \mathfrak{U}^{\beta}(n) \hbox{ and } \mathbf{P} \in \mathfrak{U}^{\beta}(m).
\end{equation}
When $m = 1$; $\mathcal{V}^{\beta}_{1,n}$ defines the unit sphere in
$\mathfrak{F}^{n}$. This is, of course, an $(n-1)\beta$- dimensional surface in
$\mathfrak{F}^{n}$. When $m = n$ and denoting $\mathbf{H}_{1}$ by $\mathbf{H}$,
$(\mathbf{H}^{*}d\mathbf{H})$ is termed the \emph{Haar measure} on
$\mathfrak{U}^{\beta}(m)$ and defines an invariant differential form of a unique
measure $\nu$ on $\mathfrak{U}^{\beta}(m)$ given by
$$
  \nu (\mathfrak{M}) = \int_{\mathfrak{M}}(\mathbf{H}^{*}d\mathbf{H}).
$$
It is unique in the sense that any other invariant measure on
$\mathfrak{U}^{\beta}(m)$ is a finite multiple of $\nu$ and invariant because is
invariant under left and right translations, that is
$$
  \nu(\mathbf{Q}\mathfrak{M}) = \nu(\mathfrak{M}\mathbf{Q}) = \nu(\mathfrak{M}),
  \quad \forall \mathbf{Q} \in \mathfrak{U}^{\beta}(m).
$$
The surface area or volume of the Stiefel manifold $\mathcal{V}^{\beta}_{m,n}$ is
\begin{equation}\label{vol}
    \Vol(\mathcal{V}^{\beta}_{m,n}) = \int_{\mathbf{H}_{1} \in
  \mathcal{V}^{\beta}_{m,n}} (\mathbf{H}^{*}_{1}d\mathbf{H}_{1}) =
  \frac{2^{m}\pi^{mn\beta/2}}{\Gamma^{\beta}_{m}[n\beta/2]},
\end{equation}
and therefore
$$
  (d\mathbf{H}_{1}) = \frac{1}{\Vol\left(\mathcal{V}^{\beta}_{m,n}\right)}
    (\mathbf{H}_{1}^{*}d\mathbf{H}_{1}) = \frac{\Gamma^{\beta}_{m}[n\beta/2]}{2^{m}
    \pi^{mn\beta/2}}(\mathbf{H}_{1}^{*}d\mathbf{H}_{1}).
$$
is the \emph{normalised invariant measure on} $\mathcal{V}^{\beta}_{m,n}$ and
$(d\mathbf{H})$, i.e. with  $(m = n)$, it defines the \emph{normalised Haar measure}
on $\mathfrak{U}^{\beta}(m)$. In (\ref{vol}), $\Gamma^{\beta}_{m}[a]$ denotes the
\emph{multivariate Gamma function} for the space $\mathfrak{S}_{m}^{\beta}$, and is
defined by
\begin{eqnarray}
  \Gamma_{m}^{\beta}[a] &=& \displaystyle\int_{\mathbf{A} \in \mathfrak{P}_{m}^{\beta}}
  \etr\{-\mathbf{A}\} |\mathbf{A}|^{a-(m-1)\beta/2 - 1}(d\mathbf{A}) \nonumber\\
    &=& \pi^{m(m-1)\beta/4}\displaystyle\prod_{i=1}^{m} \Gamma[a-(i-1)\beta/2]\nonumber\\
    &=& \pi^{m(m-1)\beta/4}\displaystyle\prod_{i=1}^{m} \Gamma[a-(m-i)\beta/2],
\end{eqnarray}
where $\etr\{\cdot\} = \exp\{\tr(\cdot)\}$, $|\cdot|$ denotes the determinant and
$\re(a)
> (m-1)\beta/2$, see \citet{gr:87}. This can be obtained as a particular case of the
\emph{generalised gamma function of weight $\kappa$} for the space
$\mathfrak{S}^{\beta}_{m}$ with $\kappa = (k_{1}, k_{2}, \dots, k_{m})$, $k_{1}\geq
k_{2}\geq \cdots \geq k_{m} \geq 0$, taking $\kappa =(0,0,\dots,0)$ and which for
$\re(a) \geq (m-1)\beta/2 - k_{m}$ is defined by, see \citet{gr:87},
\begin{eqnarray}
  \Gamma_{m}^{\beta}[a,\kappa] &=& \displaystyle\int_{\mathbf{A} \in \mathfrak{P}_{m}^{\beta}}
  \etr\{-\mathbf{A}\} |\mathbf{A}|^{a-(m-1)\beta/2 - 1} q_{\kappa}(\mathbf{A}) (d\mathbf{A}) \nonumber\\
    &=& \pi^{m(m-1)\beta/4}\displaystyle\prod_{i=1}^{m} \Gamma[a + k_{i}
    -(i-1)\beta/2]\nonumber\\
    &=& \pi^{m(m-1)\beta/4}\displaystyle\prod_{i=1}^{m} \Gamma[a + k_{i}
    -(m-i)\beta/2]\nonumber\\ \label{gammagen1}
    &=& [a]_{\kappa}^{\beta} \Gamma_{m}^{\beta}[a],
\end{eqnarray}
where for $\mathbf{A} \in \mathfrak{S}_{m}^{\beta}$
\begin{equation}\label{hwv}
    q_{\kappa}(\mathbf{A}) = |\mathbf{A}_{m}|^{k_{m}}\prod_{i = 1}^{m-1}|\mathbf{A}_{i}|^{k_{i}-k_{i+1}}
\end{equation}
with $\mathbf{A}_{p} = (a_{rs})$, $r,s = 1, 2, \dots, p$, $p = 1,2, \dots, m$ is
termed the \emph{highest weight vector}, see \citet{gr:87}.
\begin{remark}
Let $\mathcal{P}(\mathfrak{S}_{m}^{\beta})$ denote the algebra of all polynomial
functions on $\mathfrak{S}_{m}^{\beta}$, and
$\mathcal{P}_{k}(\mathfrak{S}_{m}^{\beta})$ the subspace of homogeneous polynomials
of degree $k$ and let $\mathcal{P}^{\kappa}(\mathfrak{S}_{m}^{\beta})$ be an
irreducible subspace of $\mathcal{P}(\mathfrak{S}_{m}^{\beta})$ such that
$$
  \mathcal{P}_{k}(\mathfrak{S}_{m}^{\beta}) = \sum_{\kappa}\bigoplus
  \mathcal{P}^{\kappa}(\mathfrak{S}_{m}^{\beta}).
$$
Note that $q_{\kappa}$ is a homogeneous polynomial of degree $k$, moreover
$q_{\kappa} \in \mathcal{P}^{\kappa}(\mathfrak{S}_{m}^{\beta})$, see \citet{gr:87}.
\end{remark}
In (\ref{gammagen1}), $[a]_{\kappa}^{\beta}$ denotes the generalised Pochhammer
symbol of weight $\kappa$, defined as
$$
  [a]_{\kappa}^{\beta} = \prod_{i = 1}^{m}(a-(i-1)\beta/2)_{k_{i}} = \frac{\pi^{m(m-1)\beta/4}
    \displaystyle\prod_{i=1}^{m} \Gamma[a + k_{i} -(i-1)\beta/2]}{\Gamma_{m}^{\beta}[a]} =
  \frac{\Gamma_{m}^{\beta}[a,\kappa]}{\Gamma_{m}^{\beta}[a]},
$$
where $\re(a) > (m-1)\beta/2 - k_{m}$ and
$$
  (a)_{i} = a (a+1)\cdots(a+i-1),
$$
is the standard Pochhammer symbol.

A variant of the generalised gamma function of weight $\kappa$ is obtained from
\citet{k:66} and is defined as%
\begin{eqnarray}
  \Gamma_{m}^{\beta}[a,-\kappa] &=& \displaystyle\int_{\mathbf{A} \in \mathfrak{P}_{m}^{\beta}}
    \etr\{-\mathbf{A}\} |\mathbf{A}|^{a-(m-1)\beta/2 - 1} q_{\kappa}(\mathbf{A}^{-1})
    (d\mathbf{A}) \nonumber\\
    &=& \pi^{m(m-1)\beta/4}\displaystyle\prod_{i=1}^{m} \Gamma[a - k_{i}
    -(m-i)\beta/2] \nonumber\\
    &=& \pi^{m(m-1)\beta/4}\displaystyle\prod_{i=1}^{m} \Gamma[a - k_{i}
    -(i-1)\beta/2] \nonumber \\ \label{gammagen2}
    &=& \displaystyle\frac{(-1)^{k} \Gamma_{m}^{\beta}[a]}{[-a +(m-1)\beta/2
    +1]_{\kappa}^{\beta}} ,
\end{eqnarray}
where $\re(a) > (m-1)\beta/2 + k_{1}$.

The two expressions of $\Gamma_{m}^{\beta}[a,]$, $\Gamma_{m}^{\beta}[a,\kappa]$ and
$\Gamma_{m}^{\beta}[a,-\kappa]$ as the product of ordinary gamma functions are
obtained using the proofs corresponding to $\mathbf{A} = \mathbf{TT}^{*}$ and
$\mathbf{A} = \mathbf{T}^{*}\mathbf{T}$ with the corresponding Jacobian given in
Lemma \ref{lemch}. Alternatively, note that for any function $g(y)$
\begin{equation}\label{prod1}
    \prod_{i = 1}^{q}g(x + i - 1) = \prod_{i = 1}^{q}g(x + q - i),
\end{equation}
and
\begin{equation}\label{prod2}
    \prod_{i = 1}^{q}g(x-i+1) = \prod_{i = 1}^{q}g(x-q +i).
\end{equation}

Similarly, from \citet[p. 480]{h:55} the \emph{multivariate beta function} for the
space $\mathfrak{S}^{\beta}_{m}$, can be defined as
\begin{eqnarray}
    \mathcal{B}_{m}^{\beta}[b,a] &=& \int_{\mathbf{0}<\mathbf{S}<\mathbf{I}_{m}}
    |\mathbf{S}|^{a-(m-1)\beta/2-1} |\mathbf{I}_{m} - \mathbf{S}|^{b-(m+1)\beta/2-1}
    (d\mathbf{S}) \nonumber\\
    &=& \int_{\mathbf{R} \in \mathfrak{P}_{m}^{\beta}} |\mathbf{R}|^{a-(m-1)\beta/2-1}
    |\mathbf{I}_{m} + \mathbf{R}|^{-(a+b)} (d\mathbf{R}) \nonumber\\ \label{beta}
    &=& \frac{\Gamma_{m}^{\beta}[a] \Gamma_{m}^{\beta}[b]}{\Gamma_{m}^{\beta}[a+b]},
\end{eqnarray}
where $\mathbf{R} = (\mathbf{I}-\mathbf{S})^{-1} -\mathbf{I}$, Re$(a) > (m-1)\beta/2$
and Re$(b)> (m-1)\beta/2$.

Some Jacobians in the quaternionic case are obtained in \citet{lx:09}. We now cite
some Jacobians in terms of the parameter $\beta$, based on the work of \citet{d:02}.
We also include a parameter count (or number of functionally independent variables,
\#fiv), that is, if $\mathbf{A}$ is factorised as $\mathbf{A} = \mathbf{BC}$, then
the parameter count is written as \#fiv in \textbf{A} = [\#fiv in \textbf{B}] +
[\#fiv in \textbf{C}], see \citet{d:02}.

\begin{lemma}\label{lemhlt}
Let $\mathbf{X}$ and $\mathbf{Y} \in \mathfrak{P}_{m}^{\beta}$, and let $\mathbf{Y} =
\mathbf{AXA^{*}} + \mathbf{C}$, where $\mathbf{A}$ and $\mathbf{C} \in {\mathcal
L}_{m,m}^{\beta}$ are constant matrices. Then
\begin{equation}\label{hlt}
    (d\mathbf{Y}) = |\mathbf{A}^{*}\mathbf{A}|^{\beta(m-1)/2+1} (d\mathbf{X}).
\end{equation}
\end{lemma}

\begin{lemma} [Cholesky's decomposition]\label{lemch}
Let $\mathbf{S} \in \mathfrak{P}_{m}^{\beta}$ and $\mathbf{T} \in
\mathfrak{T}_{U}^{\beta}(m)$ with $t_{ii} > 0$, $i = 1, 2, \ldots , m$. Then

$\bullet$ parameter count: $\beta m(m-1)/2 + m = \beta m(m-1)/2 + m$ and
\begin{equation}\label{ch}
    (d\mathbf{S}) = \left\{
                      \begin{array}{ll}
                        2^{m} \displaystyle\prod_{i = 1}^{m} t_{ii}^{\beta(m - i) + 1} (d\mathbf{T})
                            & \hbox{if } \ \mathbf{S} = \mathbf{T}^{*}\mathbf{T};  \\
                        2^{m} \displaystyle\prod_{i = 1}^{m} t_{ii}^{\beta(i - i) + 1} (d\mathbf{T})
                            &\hbox{if } \ \mathbf{S} = \mathbf{TT}^{*}.
                      \end{array}
                    \right.
\end{equation}
\end{lemma}

\begin{lemma}[ Spectral decomposition]\label{lemsd}
Let $\mathbf{S} \in \mathfrak{P}_{m}^{\beta}$. Then, the spectral decomposition can
be written as $\mathbf{S} = \mathbf{W}\mathbf{\Lambda W}^{*}$, where $\mathbf{W} \in
\mathfrak{U}^{\beta}(m)$ and $\mathbf{\Lambda} = \diag(\lambda_{1}, \dots,
\lambda_{m}) \in \mathfrak{D}_{m}^{1}$, with $\lambda_{1}> \cdots> \lambda_{m}>0$.
Then

$\bullet$ parameter count: $\beta m(m-1)/2 + m = [\beta m(m+1)/2 - m -(\beta -1)m] +
[m]$ and
\begin{equation}\label{sd}
    (d\mathbf{S}) = 2^{-m} \pi^{\varrho} \prod_{i < j}^{m} (\lambda_{i} - \lambda_{j})^{\beta}
    (d\mathbf{\Lambda})(\mathbf{W}^{*}d\mathbf{W}),
\end{equation}
where
$$
  \varrho = \left\{
             \begin{array}{rl}
               0, & \beta = 1; \\
               -m, & \beta = 2; \\
               -2m, & \beta = 4; \\
               -4m, & \beta = 8.
             \end{array}
           \right.
$$
\end{lemma}

\section{Integral properties of Jack polynomials}\label{sec3}

In this section we review and study diverse integral properties of Jack polynomials
for normed division algebras. However, let us first consider the following remarks
and definitions.

\begin{remark}
Note that Jack polynomials and hypergeometric functions with one or two matrix
arguments are valid for $\beta > 0$ (\citet{ke:06}), but \emph{in our case} $\beta$
denotes the \emph{real dimension} of $\mathfrak{F}$. Also, we use the parameter
$\beta$ indeed of $\alpha$ in the definition of the Jack polynomials and
hypergeometric functions, with the equivalence shown in Table \ref{table1}.
\end{remark}

Then:

Let us characterise the \emph{Jack symmetric function}
$J_{\kappa}^{(\beta)}(\lambda_{1},\ldots,\lambda_{m})$ of parameter $\beta$, see
\citet{S97}. A decreasing sequence of nonnegative integers
$\kappa=(k_{1},k_{2},\ldots)$ with only finitely many nonzero terms is said to be a
partition of $k=\sum k_{i}$. Let $\kappa$ and $\tau=(t_{1},t_{2},\ldots)$ be two
partitions of $k$. We write $\tau\leq\kappa$ if
$\sum_{i=1}^{t}t_{i}\leq\sum_{i=1}^{t}k_{i}$ for each $t$. The conjugate of $\kappa$
is $\kappa'=(k_{1}',k_{2}',\ldots)$ where $k_{i}'= \mbox{card}\{j:k_{j}\geq i\}$. The
length of $\kappa$ is $l(k)=\mbox{max}\{i:k_{i}\neq 0\}=k_{1}'$.  If $l(\kappa)\leq
m$, it is often written that $\kappa=(k_{1},k_{2},\ldots,k_{m})$.

The \emph{monomial symmetric function} $M_{\kappa}(\cdot)$ indexed by a partition
$\kappa$ can be regarded as a function of an arbitrary number of variables such that
all but a finite number are equal to $0$:  if $\lambda_{i}=0$ for $i>m\geq l(\kappa)$
then $M_{\kappa}(\lambda_{1},\ldots,\lambda_{m})=\sum \ \lambda_{1}^{\delta_{1}} \
\cdots  \ \lambda_{m}^{\delta_{m}}$, where the sum is over all distinct permutations
$\{\delta_{1},\ldots,\delta_{m}\}$ of $\{k_{1},\ldots,k_{m}\}$, and if $l(\kappa)>m$
then $M_{\kappa}(\lambda_{1},\ldots,\lambda_{m})=0$.  A symmetric function $f$ is a
linear combination of monomial symmetric functions.  If $f$ is a symmetric function
then $f(\lambda_{1},\ldots,\lambda_{m},0)=f(\lambda_{1},\ldots,\lambda_{m})$.  For
each $m\geq 1$, $f(\lambda_{1},\ldots,\lambda_{m})$ is a symmetric polynomial in $m$
variables.

Then the \emph{Jack symmetric function}
$J_{\kappa}^{(\beta)}(\lambda_{1},\ldots,\lambda_{m})$ with a parameter $\beta$,
satisfies the following conditions:
\begin{eqnarray}
    J_{\kappa}^{(\beta)}(\lambda_{1},\ldots,\lambda_{m})&=&\sum_{\tau\leq\kappa} \nu_{\kappa,\tau}(\beta)
    M_{\tau}(\lambda_{1},\ldots,\lambda_{m}),\label{eq:Condition1}\\
    J_{\kappa}^{(\beta)}(1,\ldots,1) &=& \left (\frac{2}{\beta}\right)^{k}\prod_{i=1}^{m}
    \left((m-i+1)\beta/2\right)_{k_{i}}, \label{eq:Condition2}\\
    \mathcal{D}_{2}^{\beta}J_{\kappa}^{(\beta)}(\lambda_{1},\ldots,\lambda_{m}) &=&
    \sum_{i=1}^{m}k_{i}(k_{i}-1+ \beta (m-i))J_{\kappa}^{(\beta)}(\lambda_{1},\ldots,\lambda_{m})
    \label{eq:Condition3}.
\end{eqnarray}
where
$$
  \mathcal{D}_{2}^{\beta} = \sum_{i =1}^{m}\lambda_{i}^{2}\frac{\partial^{2}}{\partial \lambda_{i}^{2}}
    + \beta \sum_{i=1}^{m}\lambda_{i}^{2}\sum_{j\neq i}\frac{1}{\lambda_{i}-\lambda_{j}}
    \frac{\partial}{\partial \lambda_{i}}.
$$

Here, the constants $\nu_{\kappa,\tau}(\beta)$ do not dependent on $\lambda_{i}'s$
but on $\kappa$ and $\tau$. Note that if $m<l(\kappa)$ then
$J_{\kappa}^{(\beta)}(\lambda_{1},\ldots,\lambda_{m})=0$. The conditions include the
case $\beta =0$ and then $J_{\kappa}^{(0)}(\lambda_{1},\ldots,\lambda_{m})=e_{\kappa
'}\prod_{i=1}^{m}(m-i+1)^{k_{i}}$, where
$e_{\kappa}(\lambda_{1},\ldots,\lambda_{m})=\prod_{i=1}^{l(\kappa)}e_{k_{i}}(\lambda_{1},\ldots,\lambda_{m})$
are the elementary symmetric functions indexed by partitions $\kappa$, if $m\geq
l(\kappa)$ then
$e_{r}(\lambda_{1},\ldots,\lambda_{m})=\sum_{i_{1}<i_{2}<\cdots<i_{r}}\lambda_{i_{1}}\cdots
\lambda_{i_{r}}$, and if $m< l(\kappa)$ then
$e_{r}(\lambda_{1},\ldots,\lambda_{m})=0$, see \citet{S97}.

Now, from \citet{ke:06}, the Jack functions
$$
  J_{\kappa}^{(\beta)}(\mathbf{X}) = J_{\kappa}^{(\beta)}(\lambda_{1},\ldots,\lambda_{m}),
$$
wher $\lambda_{1},\ldots,\lambda_{m}$ are the eigenvalues of the matrix $\mathbf{X}
\in \mathfrak{S}_{m}^{\beta}$, can be normalised in such a way that
\begin{equation}\label{jptr}
    \sum_{\kappa}C_{\kappa}^{\beta}(\mathbf{X}) = (\tr(\mathbf{X}))^{k},
\end{equation}
or equivalently, such that
\begin{equation}\label{jpexp}
    \sum_{k = 1}^{\infty}\sum_{\kappa}\frac{C_{\kappa}^{\beta}(\mathbf{X})}{k!} = \etr\{\mathbf{X}\},
\end{equation}
where $C_{\kappa}^{\beta}(\mathbf{X})$ denotes the \emph{Jack polynomials} (for
simplicity, we have replaced $(\beta)$ by $\beta$ as the superindex for the Jack
polynomials). These are related to the Jack functions by
\begin{equation}\label{eq1}
    C_{\kappa}^{\beta}(\mathbf{X})= \frac{2^{k} k!}{\beta^{k} \nu_{\kappa}}
        J_{\kappa}^{(\beta)}(\mathbf{X}),
\end{equation}
where
$$
  \nu_{\kappa} = \prod_{(i,j) \in \kappa} h_{*}^{\kappa}(i,j) h^{*}_{\kappa}(i,j),
$$
and $h_{*}^{\kappa}(i,j) = k_{j}-i + 2(k_{i}-j+1)/\beta$ and $h^{*}_{\kappa}(i,j) =
k_{j}-i+1 + 2(k_{i}-j)/\beta$ are the upper and lower hook lengths at $(i,j) \in
\kappa$, respectively. Also, observe that for $\mathbf{X} = \mathbf{S}^{*}
\mathbf{S}$ and $\mathbf{Y}= \mathbf{W}^{*} \mathbf{W}$ we have
\begin{equation}
    C_{\kappa}^{\beta}(\mathbf{W}\mathbf{XW}^{*}) =
     C_{\kappa}^{\beta}(\mathbf{S}\mathbf{YS}^{*}).
\end{equation}
In particular for $\mathbf{A}^{1/2}$ such that $\left(\mathbf{A}^{1/2}\right)^{2} =
\mathbf{A}$
\begin{equation}
    C_{\kappa}^{\beta}(\mathbf{Y}^{1/2}\mathbf{XY}^{1/2}) =
        C_{\kappa}^{\beta}(\mathbf{X}^{1/2}\mathbf{YX}^{1/2}).
\end{equation}
Therefore, given that $\mathbf{XY}$, $\mathbf{YX}$,
$\mathbf{Y}^{1/2}\mathbf{XY}^{1/2}$ and $\mathbf{X}^{1/2}\mathbf{YX}^{1/2}$ all have
the same eigenvalues, we opt for convenience of notation rather than strict adherence
to rigor, and write $C_{\kappa}^{\beta}(\mathbf{XY})$ or
$C_{\kappa}^{\beta}(\mathbf{YX})$ rather than
$C_{\kappa}^{\beta}(\mathbf{Y}^{1/2}\mathbf{XY}^{1/2})$, even though $\mathbf{XY}$ or
$\mathbf{YX}$ need not lie in $\mathfrak{S}_{m}^{\beta}$. Note that
\begin{equation}\label{jpeq4}
    C_{\kappa}^{\beta}(\mathbf{Z}^{1/2}\mathbf{XZ}^{1/2}\mathbf{Y}) =
        C_{\kappa}^{\beta}(\mathbf{XZ}^{1/2}\mathbf{Y}\mathbf{Z}^{1/2}),
\end{equation}
for all $\mathbf{X}, \mathbf{Y} \in \mathfrak{S}_{m}^{\beta}$ and $\mathbf{Z} \in
\mathfrak{P}_{m}^{\beta}$. From \citet[Equation 4.8(2) and Definition 5.3]{gr:87} we
have
\begin{equation}\label{jpq}
    C_{\kappa}^{\beta}(\mathbf{X}) = C_{\kappa}^{\beta}(\mathbf{I})
    \int_{\mathbf{H} \in \mathfrak{U}^{\beta}(m)} q_{\kappa}(\mathbf{H}^{*}\mathbf{XH})(d\mathbf{H})
\end{equation}
for all $\mathbf{X} \in \mathfrak{S}_{m}^{\beta}$; where $(d\mathbf{H})$ is the
normalised Haar measure on $\mathfrak{U}^{\beta}(m)$.  Finally,  for the $c$ constant
we have that $C_{\kappa}^{\beta}(c\mathbf{X}) = c^{k}C_{\kappa}^{\beta}(\mathbf{X})$.

Some basic integral properties are cited below. For this purpose, we utilise the
complexification $\mathfrak{S}_{m}^{\beta, \mathfrak{C}} = \mathfrak{S}_{m}^{\beta} +
i \mathfrak{S}_{m}^{\beta}$ of $\mathfrak{S}_{m}^{\beta}$. That is,
$\mathfrak{S}_{m}^{\beta, \mathfrak{C}}$ consist of all matrices $\mathbf{X} \in
(\mathfrak{F^{\mathfrak{C}}})^{m \times m}$ of the form $\mathbf{Z} = \mathbf{X} +
i\mathbf{Y}$, with $\mathbf{X}, \mathbf{Y} \in \mathfrak{S}_{m}^{\beta}$. We refer to
$\mathbf{X} = \re(\mathbf{Z})$ and $\mathbf{Y} = \im(\mathbf{Z})$ as the \emph{real
and imaginary parts} of $\mathbf{Z}$, respectively. The \emph{generalised right
half-plane} $\mathbf{\Phi} = \mathfrak{P}_{m}^{\beta} + i \mathfrak{S}_{m}^{\beta}$
in $\mathfrak{S}_{m}^{\beta,\mathfrak{C}}$ consists of all $\mathbf{Z} \in
\mathfrak{S}_{m}^{\beta,\mathfrak{C}}$ such that $\re(\mathbf{Z}) \in
\mathfrak{P}_{m}^{\beta}$, see \citet[p. 801]{gr:87}.

For any $\mathbf{X}, \mathbf{Y} \in \mathfrak{S}_{m}^{\beta}$,
\begin{equation}\label{jpeq1}
    \int_{\mathbf{H} \in \mathfrak{U}^{\beta}(m)} C_{\kappa}^{\beta}(\mathbf{XH}^{*}\mathbf{YH})
    (d\mathbf{H}) = \frac{C_{\kappa}^{\beta}(\mathbf{X})C_{\kappa}^{\beta}(\mathbf{Y})}
    {C_{\kappa}^{\beta}(\mathbf{I})}.
\end{equation}

For all $\mathbf{R} \in \mathfrak{S}_{m}^{\beta}$, $\mathbf{Z} \in \mathbf{\Phi}$ and
$\re(a) > (m-1)\beta/2 -k_{m}$,
\begin{eqnarray}
  \int_{\mathbf{X} \in \mathfrak{P}_{m}^{\beta}} \etr\{-\mathbf{XZ}\} |\mathbf{X}|^{a-(m-1)\beta/2-1}
    C_{\kappa}^{\beta}(\mathbf{XR})(d\mathbf{X}) \hspace{-3.5cm} \nonumber\\ \label{jpeq2}
    &=& \Gamma_{m}^{\beta}[a, \kappa]
    |\mathbf{Z}|^{-a} C_{\kappa}^{\beta}(\mathbf{RZ}^{-1})\nonumber\\
    &=& [a]_{\kappa}^{\beta} \Gamma_{m}^{\beta}[a]
    |\mathbf{Z}|^{-a} C_{\kappa}^{\beta}(\mathbf{RZ}^{-1}).
\end{eqnarray}
\begin{remark}
   In general, the result (\ref{jpeq2}) has been established under the condition, $\re(a) >
   (m-1)\beta/2$, see \citet{c:63}, \citet{m:82}, \citet{rva:05b} and \citet{lx:09}, but
   in reality the correct condition is $\re(a) > (m-1)\beta/2 -k_{m}$. This fact is
   immediate, observing that $[a]_{\kappa}^{\beta} \Gamma_{m}^{\beta}[a] = \Gamma_{m}^{\beta}[a,
   \kappa]$ and the different expressions for $\Gamma_{m}^{\beta}[a, \kappa]$ in (\ref{gammagen1}).
\end{remark}
Let $\re(a) > (m-1)\beta/2$ and $\re(b) > (m-1)\beta/2$. Then
\begin{eqnarray}
  \int_{\mathbf{0} < \mathbf{X} < \mathbf{I}} |\mathbf{X}|^{a-(m-1)\beta/2-1} |\mathbf{I} -
    \mathbf{X}|^{b-(m-1)\beta/2-1} C_{\kappa}^{\beta}(\mathbf{XR})(d\mathbf{X})
    \hspace{-4cm} \nonumber\\\label{jpeq3}
    &=& \frac{\Gamma_{m}^{\beta}[a, \kappa]\Gamma_{m}^{\beta}[b]}{\Gamma_{m}^{\beta}[a + b, \kappa]}
    C_{\kappa}^{\beta}(\mathbf{R}) \nonumber\\
    &=& \frac{[a]_{\kappa}^{\beta}\mathcal{B}_{m}^{\beta}[a,b]}{[a + b]_{m}^{\beta}}
    C_{\kappa}^{\beta}(\mathbf{R}),
\end{eqnarray}
for all $\mathbf{R} \in \mathfrak{S}_{m}^{\beta,\mathfrak{C}}$; see \citet[Theorems
5.5 and 5.9 and Corollary 5.10]{gr:87} for real, complex and quaternion cases.
\begin{remark}
   Observe that result (\ref{jpeq3}) was established under the conditions $\re(a) >
   (m-1)\beta/2$ and $\re(b) > (m-1)\beta/2$, see \citet{c:63}, \citet{m:82}, \citet{rva:05b} and \citet{lx:09}, but
   the correct conditions are in fact $\re(a) > (m-1)\beta/2 - k_{m}$ and $\re(b) > (m-1)\beta/2$. This fact is
   verified by observing that $[a]_{\kappa}^{\beta} \Gamma_{m}^{\beta}[a] = \Gamma_{m}^{\beta}[a,
   \kappa]$ and the different expressions for $\Gamma_{m}^{\beta}[a, \kappa]$ in (\ref{gammagen1}).
\end{remark}

We now extend several integral properties of zonal polynomials in the real and
complex cases to normed division algebras. Our first result is a generalisation of
one studied by \citet{tfd:89} for real case, see also \citet{cldggf:09}. From this
result, we can obtain diverse particular integral properties of Jack polynomials.

\begin{theorem}\label{teo1}
Let $\mathbf{Z} \in \mathbf{\Phi}$ and $\mathbf{U} \in \mathfrak{S}_{m}^{\beta}$.
Assume $\gamma = \int_{0}^{\infty}f(z) z^{am-k-1} dz < \infty$. Then
\begin{eqnarray}\label{runze21}
    \int_{\mathbf{X} \in \mathfrak{P}_{m}^{\beta}} f(\tr \mathbf{XZ}) |\mathbf{X}|^{a -(m-1)
    \beta/2-1} C_{\kappa}^{\beta}\left(\mathbf{X}^{-1}\mathbf{U} \right)
    (d\mathbf{X}) \hspace{3cm} \nonumber\\\label{jpe1}
    =\displaystyle\frac{\Gamma_{m}^{\beta}[a, - \kappa]}{\Gamma[am-k]}
     |\mathbf{Z}|^{-a} C_{\kappa}^{\beta}(\mathbf{UZ})\cdot \gamma,
\end{eqnarray}
for $\re(a) > (m-1)\beta/2 + k_{1}$, and
\begin{eqnarray}\label{runze22}
    \int_{\mathbf{X} \in \mathfrak{P}_{m}^{\beta}} f(\tr \mathbf{XZ}) |\mathbf{X}|^{a -(m-1)
    \beta/2-1} C_{\kappa}^{\beta}\left(\mathbf{X}\mathbf{U} \right)
    (d\mathbf{X}) \hspace{3cm} \nonumber\\\label{jpe2}
    =\displaystyle\frac{\Gamma_{m}^{\beta}[a, \kappa]}{\Gamma[am+k]}
     |\mathbf{Z}|^{-a} C_{\kappa}^{\beta}(\mathbf{UZ}^{-1})\cdot \vartheta,
\end{eqnarray}
where $\vartheta = \int_{0}^{\infty}f(z) z^{am+k-1} dz < \infty$,  $\re(a)
> (m-1)\beta/2 - k_{m}$ and $\kappa = (k_{1}, \dots, k_{m})$ and $k =
k_{1} + \cdots + k_{m}$.
\end{theorem}
\begin{proof}
Denote the left side of (\ref{runze21}) by $I(\mathbf{U},\mathbf{Z})$. By (\ref{jpq})
and interchange of order on integration
\begin{eqnarray*}
  I(\mathbf{I},\mathbf{I}) &=& \int_{\mathbf{X} \in \mathfrak{P}_{m}^{\beta}} f(\tr \mathbf{X})
    |\mathbf{X}|^{a -(m-1) \beta/2-1} C_{\kappa}^{\beta}\left(\mathbf{X}^{-1} \right)
    (d\mathbf{X}) \\
    &=&  C_{\kappa}^{\beta}(\mathbf{I}) \int_{\mathbf{X} \in \mathfrak{P}_{m}^{\beta}} f(\tr \mathbf{X})
    |\mathbf{X}|^{a -(m-1) \beta/2-1}\\
    & & \hspace{3.5cm} \times \left(\int_{\mathbf{H} \in \mathfrak{U}^{\beta}(m)}
    q_{\kappa}\left(\mathbf{H}^{*}\mathbf{X}^{-1}\mathbf{H} \right)(d\mathbf{H})\right) (d\mathbf{X})\\
    &=& C_{\kappa}^{\beta}(\mathbf{I}) \int_{\mathbf{X} \in \mathfrak{P}_{m}^{\beta}} f(\tr \mathbf{X})
    |\mathbf{X}|^{a -(m-1) \beta/2-1}
    q_{\kappa}\left(\mathbf{X}^{-1}\right)(d\mathbf{X}).
\end{eqnarray*}
Let $\mathbf{X} = \mathbf{TT}^{*}$, from Lemma \ref{lemch}
$$
  (d\mathbf{X})= 2^{m} \prod_{i = 1}^{m} t_{ii}^{\beta(i - i) + 1} (d\mathbf{T}).
$$
Then
$$
  \mathcal{I}(\mathbf{I},\mathbf{I}) = 2^{m} C_{\kappa}^{\beta}(\mathbf{I}) \build{\int \cdots \int}{}
    {\build{}{0 < t_{ii} < \infty} {-\infty < t_{ij} < \infty}}{}f\left(\sum_{i = 1}^{m}
    t_{ij}^{2}\right)\prod_{i=1}^{m}(t_{ii})^{2(a - k_{i} - (m -i)\beta/2)-1}
    (d\mathbf{T}).
$$
Applying the \citet[Lemma 2.4.3, p. 51]{fz:90} we obtain
\begin{eqnarray*}
  \mathcal{I}(\mathbf{I},\mathbf{I}) &=& C_{\kappa}^{\beta}(\mathbf{I})
  \frac{\pi^{m(m-1)\beta/4}\displaystyle\prod_{i = 1}^{m}\Gamma[a-k_{i}-(m-i)\beta/2]}{\Gamma[am - k]} \cdot
  \gamma \\
  &=& C_{\kappa}^{\beta}(\mathbf{I}) \frac{\Gamma_{m}^{\beta}[a,-\kappa]}{\Gamma[am - k]} \cdot \gamma,
\end{eqnarray*}
with $\gamma = \int_{0}^{\infty}f(z) z^{am-k-1} dz < \infty$.

Next, since the function $\mathcal{I}(\mathbf{U},\mathbf{I})$ is invariant under
$\mathfrak{U}^{\beta}(m)$ and in $\mathcal{P}^{\kappa}(\mathfrak{S}_{m}^{\beta})$,
there exists a constant $d$ such that $\mathcal{I}(\mathbf{U},\mathbf{I}) = d
C_{\kappa}^{\beta}(\mathbf{U})$, $\mathbf{U} \in \mathfrak{S}_{m}^{\beta}$. It is
obvious that $d = \mathcal{I}(\mathbf{I},\mathbf{I})/C_{\kappa}^{\beta}(\mathbf{I})$,
then
$$
  \mathcal{I}(\mathbf{U},\mathbf{I}) =  \frac{\Gamma_{m}^{\beta}[a,-\kappa]}{\Gamma[am - k]}
  C_{\kappa}^{\beta}(\mathbf{U})\cdot \gamma.
$$
Now, let $\mathbf{Z} \in \mathfrak{P}_{m}^{\beta}$ and make the change of variable
$\mathbf{X} \rightarrow \mathbf{Z}^{-1/2}\mathbf{XZ}^{-1/2}$ in the integral defining
$\mathcal{I}(\mathbf{U},\mathbf{Z})$. Then by (\ref{jpeq4})
$$
  \mathcal{I}(\mathbf{U},\mathbf{Z}) = |\mathbf{Z}|^{-a}
  \mathcal{I}(\mathbf{Z}^{1/2} \mathbf{U}\mathbf{Z}^{1/2},\mathbf{I}),
$$
and hence
$$
  \mathcal{I}(\mathbf{U},\mathbf{Z}) = \frac{\Gamma_{m}^{\beta}[a,-\kappa]}{\Gamma[am - k]}
  |\mathbf{Z}|^{-a}C_{\kappa}^{\beta}(\mathbf{Z}^{1/2}\mathbf{U}\mathbf{Z}^{1/2})\cdot
  \gamma.
$$
Therefore, for $\mathbf{Z} \in \mathfrak{P}_{m}^{\beta}$ and $\mathbf{U} \in
\mathfrak{S}_{m}^{\beta}$
$$
  \mathcal{I}(\mathbf{U},\mathbf{Z}) = \frac{\Gamma_{m}^{\beta}[a,-\kappa]}{\Gamma[am - k]}
  |\mathbf{Z}|^{-a}C_{\kappa}^{\beta}(\mathbf{U}\mathbf{Z})\cdot
  \gamma.
$$
The result in (\ref{jpe1}) now follows by analytic continuation in $\mathbf{Z}$ from
$\mathfrak{P}_{m}^{\beta}$ to $\mathbf{\Phi} = \mathfrak{P}_{m}^{\beta} +
i\mathfrak{S}_{m}^{\beta}$. The result in (\ref{jpe2}) is obtained in a similar way.
\end{proof}

Now, by definition if $\kappa = 0$ then $[a]_{\kappa}^{\beta} = 1$ and
$C_{\kappa}^{\beta}(\mathbf{X}) = 1$ from where:

\begin{corollary}\label{coro1}
Let $\mathbf{Z} \in \mathbf{\Phi} \in \mathfrak{S}_{m}^{\beta}$. Assume $\gamma =
\int_{0}^{\infty}f(z) z^{am-1} dz < \infty$. Then
\begin{equation}
    \int_{\mathbf{X} \in \mathfrak{P}_{m}^{\beta}} f(\tr \mathbf{XZ}) |\mathbf{X}|^{a -(m-1)
    \beta/2-1} (d\mathbf{X}) = \frac{\Gamma_{m}^{\beta}[a]}{\Gamma[am]}
     |\mathbf{Z}|^{-a} \cdot \gamma,
\end{equation}
for $\re(a) > (m-1)\beta/2$.
\end{corollary}

If we take $f(y) = \exp\{- y\}$ in Corollary \ref{coro1} we obtain
\begin{equation}\label{gengamma}
    \int_{\mathbf{X} \in \mathfrak{P}_{m}^{\beta}} \etr\{ -\mathbf{XZ}\} |\mathbf{X}|^{a -(m-1)
    \beta/2-1} (d\mathbf{X}) = \Gamma_{m}^{\beta}[a] |\mathbf{Z}|^{-a},
\end{equation}
and if $\mathbf{Z} = \mathbf{I}$ we obtain the multivariate gamma function for the
space $\mathfrak{S}_{m}^{\beta}$.

Other particular results of Theorem \ref{teo1} are summarised below.

\begin{corollary}\label{coro2} Let $\mathbf{Z} \in \mathbf{\Phi}$ and $\mathbf{U} \in \mathfrak{S}_{m}^{\beta}$.
\begin{eqnarray}
  \int_{\mathbf{X} \in \mathfrak{P}_{m}^{\beta}} \etr\{- \mathbf{XZ}\} |\mathbf{X}|^{a -(m-1)
    \beta/2-1} C_{\kappa}^{\beta}\left(\mathbf{X}^{-1}\mathbf{U} \right)
    (d\mathbf{X}) \hspace{3cm} \nonumber \\
    =\Gamma_{m}^{\beta}[a, - \kappa]|\mathbf{Z}|^{-a} C_{\kappa}^{\beta}(\mathbf{UZ}),
\end{eqnarray}
for $\re(a) > (m-1)\beta/2 + k_{1}$.
\end{corollary}
\begin{proof}
This is obtained by taking $f(y) = \exp \{-y\}$ in (\ref{runze21}).
\end{proof}

\begin{corollary}\label{coro3} Let $\mathbf{Z} \in \mathbf{\Phi}$ and $\mathbf{U} \in \mathfrak{S}_{m}^{\beta}$
and $j \in \Re$, such that $\re(ma + j - k) > 0$, then
\begin{eqnarray}
    \int_{\mathbf{X} \in \mathfrak{P}_{m}^{\beta}} \etr\{- \mathbf{XZ}\} (\tr \mathbf{XZ})^{j}
    |\mathbf{X}|^{a -(m-1)\beta/2-1} C_{\kappa}^{\beta}\left(\mathbf{X}^{-1}\mathbf{U} \right)
    (d\mathbf{X}) \hspace{2cm} \nonumber\\
    = \frac{\Gamma_{m}^{\beta}[a, - \kappa] \Gamma[ma + j - k]}{\Gamma[ma-k]}|\mathbf{Z}|^{-a}
    C_{\kappa}^{\beta}(\mathbf{UZ}),
\end{eqnarray}
for $\re(a) > (m-1)\beta/2 + k_{1}$. And if $j$ is such that $\re(ma + j + k) > 0$,
then
\begin{eqnarray}
    \int_{\mathbf{X} \in \mathfrak{P}_{m}^{\beta}} \etr\{- \mathbf{XZ}\} (\tr \mathbf{XZ})^{j}
    |\mathbf{X}|^{a -(m-1) \beta/2-1} C_{\kappa}^{\beta}\left(\mathbf{X}\mathbf{U} \right)
    (d\mathbf{X}) \hspace{2cm} \nonumber\\
    = \frac{\Gamma_{m}^{\beta}[a, \kappa] \Gamma[ma + j + k]}{\Gamma[ma+k]}|\mathbf{Z}|^{-a}
    C_{\kappa}^{\beta}(\mathbf{UZ}^{-1}),
\end{eqnarray}
for $\re(a) > (m-1)\beta/2 - k_{m}$.
\end{corollary}
\begin{proof}
The desired result is obtained by taking $f(y) = \exp \{-y\} y^{j}$ in Theorem
\ref{teo1}.
\end{proof}

\begin{corollary}\label{coro4} Let $\mathbf{Z} \in \mathbf{\Phi}$ and $\mathbf{U} \in \mathfrak{S}_{m}^{\beta}$
and $\eta > 0$ then
\begin{eqnarray}
    \int_{\mathbf{X} \in \mathfrak{P}_{m}^{\beta}} (1 + 2\eta^{-1}\tr \mathbf{XZ})^{-\beta(am+\eta)}
    |\mathbf{X}|^{a -(m-1)\beta/2-1} C_{\kappa}^{\beta}\left(\mathbf{X}^{-1}\mathbf{U} \right)
    (d\mathbf{X}) \hspace{0.5cm} \nonumber\\
    = \frac{\Gamma_{m}^{\beta}[a, - \kappa] \Gamma[(\beta-1)am + \beta \eta + k]}{(2\eta^{-1})^{am-k}
    \Gamma[\beta(ma+\eta)]}|\mathbf{Z}|^{-a} C_{\kappa}^{\beta}(\mathbf{UZ}),
\end{eqnarray}
for $\re(a) > (m-1)\beta/2 + k_{1}$. And
\begin{eqnarray}
    \int_{\mathbf{X} \in \mathfrak{P}_{m}^{\beta}} (1 + 2\eta^{-1}\tr \mathbf{XZ})^{-\beta(am+\eta)}
    |\mathbf{X}|^{a -(m-1)\beta/2-1} C_{\kappa}^{\beta}\left(\mathbf{X}\mathbf{U} \right)
    (d\mathbf{X}) \hspace{0.5cm} \nonumber\\
    = \frac{\Gamma_{m}^{\beta}[a, \kappa] \Gamma[(\beta-1)am + \beta \eta - k]}{(2\eta^{-1})^{am+k}
    \Gamma[\beta(ma+\eta)]}|\mathbf{Z}|^{-a} C_{\kappa}^{\beta}(\mathbf{UZ}^{-1}),
\end{eqnarray}
for $\re(a) > (m-1)\beta/2 - k_{m}$.
\end{corollary}
\begin{proof}
The desired result is obtained  by taking $f(y) = (1+2\eta^{-1} y)^{-\beta(am+\eta)}$
in Theorem \ref{teo1}.
\end{proof}

Many other interesting particular cases of Theorem \ref{teo1} can be found, for
example by defining $f(\tr \mathbf{XZ})$ as the kernel of matrix variate generalised
Wishart distributions, see \citet{fz:90} and \citet{gv:93}.

Important analogues of the beta function integral are given in the following
theorems. Theorem \ref{teo2} is discussed by \citet{k:66} in the real case.

\begin{theorem}\label{teo2}
If $\mathbf{R} \in \mathfrak{S}_{m}^{\beta,\mathfrak{C}}$, then
\begin{eqnarray}
  \int_{\mathbf{X} \in \mathfrak{P}_{m}^{\beta}} |\mathbf{X}|^{a-(m-1)\beta/2-1}
  |\mathbf{I} + \mathbf{X}|^{-(a+b)} C_{\kappa}^{\beta}(\mathbf{RX}^{-1})(d\mathbf{X})
  \hspace{2cm} \nonumber\\ \label{r1teo2}
  = \frac{\Gamma_{m}^{\beta}[a, -\kappa] \Gamma_{m}^{\beta}[b, \kappa]}{\Gamma_{m}^{\beta}[a + b]}
  C_{\kappa}^{\beta}(\mathbf{R}),
\end{eqnarray}
for $\re(a) > (m-1)\beta/2 + k_{1}$ and $\re(b) > (m-1)\beta/2 - k_{m}$. And
\begin{eqnarray}
  \int_{\mathbf{X} \in \mathfrak{P}_{m}^{\beta}} |\mathbf{X}|^{a-(m-1)\beta/2-1}
  |\mathbf{I} + \mathbf{X}|^{-(a+b)} C_{\kappa}^{\beta}(\mathbf{RX})(d\mathbf{X})
  \hspace{2cm} \nonumber\\\label{r2teo2}
  = \frac{\Gamma_{m}^{\beta}[a, \kappa] \Gamma_{m}^{\beta}[b, -\kappa]}{\Gamma_{m}^{\beta}[a + b]}
  C_{\kappa}^{\beta}(\mathbf{R}),
\end{eqnarray}
for $\re(a) > (m-1)\beta/2 - k_{m}$ and $\re(b) > (m-1)\beta/2 + k_{1}$.
\end{theorem}
\begin{proof}
By Corollary \ref{coro2}, we have for any $\mathbf{Z} \in \mathbf{\Phi}$
\begin{eqnarray}
  \int_{\mathbf{X} \in \mathfrak{P}_{m}^{\beta}} \etr\{- \mathbf{XZ}\} |\mathbf{X}|^{a -(m-1)
    \beta/2-1} C_{\kappa}^{\beta}\left(\mathbf{X}^{-1}\mathbf{R} \right)
     |\mathbf{Z}|^{a}(d\mathbf{X}) \hspace{2cm} \nonumber \\ \label{proofteo21}
    =\Gamma_{m}^{\beta}[a, - \kappa] C_{\kappa}^{\beta}(\mathbf{RZ}).
\end{eqnarray}
Multiplying both sides of (\ref{proofteo21}) by $\etr\{-\mathbf{Z}\}|\mathbf{Z}|^{b
-(m-1)\beta/2-1}$ and integrating with respect to $\mathbf{Z}$ we have
\begin{eqnarray*}
  \int_{\mathbf{X} \in \mathfrak{P}_{m}^{\beta}} \left(\int_{\mathbf{Z} \in \mathfrak{P}_{m}^{\beta}}
    \etr\{- (\mathbf{I} +\mathbf{X})\mathbf{Z}\} |\mathbf{Z}|^{a + b -(m-1)\beta/2-1}(d\mathbf{Z})\right)\hspace{2cm}\\
    \times |\mathbf{X}|^{a -(m-1)\beta/2-1} C_{\kappa}^{\beta}\left(\mathbf{X}^{-1}\mathbf{R} \right)
    (d\mathbf{X}) \nonumber \\
\end{eqnarray*}
\begin{equation}\label{proofteo22}
    \phantom{********} = \Gamma_{m}^{\beta}[a, - \kappa]  \int_{\mathbf{Z} \in \mathfrak{P}_{m}^{\beta}}
    \etr\{-\mathbf{Z}\}|\mathbf{Z}|^{b -(m-1)\beta/2-1} C_{\kappa}^{\beta}(\mathbf{RZ})(d\mathbf{Z}).
\end{equation}
The desired result in (\ref{r1teo2}) is obtained by using Corollary \ref{coro2} and
(\ref{gengamma}) in the left and right sides of (\ref{proofteo22}), respectively. The
result in (\ref{r2teo2}) is obtained similarly.
\end{proof}

\begin{corollary}\label{coro5}
If $\mathbf{R} \in \mathfrak{S}_{m}^{\beta,\mathfrak{C}}$, then
\begin{eqnarray}
  \int_{\mathbf{0} < \mathbf{X} < \mathbf{I}} |\mathbf{X}|^{a-(m-1)\beta/2-1}
  |\mathbf{I} - \mathbf{X}|^{b-(m-1)\beta/2-1} C_{\kappa}^{\beta}(\mathbf{RX}^{-1})(d\mathbf{X})
  \hspace{2cm} \nonumber\\
  = \frac{\Gamma_{m}^{\beta}[a, -\kappa] \Gamma_{m}^{\beta}[b]}{\Gamma_{m}^{\beta}[a + b, -\kappa]}
  C_{\kappa}^{\beta}(\mathbf{R}),
\end{eqnarray}
for $\re(a) > (m-1)\beta/2 + k_{1}$ and $\re(b) > (m-1)\beta/2$.
\end{corollary}
\begin{proof}
It is obtained in a similar way to that given for (\ref{jpeq3}), see \citet{gr:87}.
\end{proof}

Now, taking $b = (m-1)\beta/2+1 > (m-1)\beta/2$, from (\ref{jpeq3}) and Corollary
\ref{coro5} we have the following result.

\begin{corollary}\label{coro6}
If $\mathbf{R} \in \mathfrak{S}_{m}^{\beta,\mathfrak{C}}$, then
\begin{eqnarray}
  \int_{\mathbf{0} < \mathbf{X} < \mathbf{I}} |\mathbf{X}|^{a-(m-1)\beta/2-1}
  C_{\kappa}^{\beta}(\mathbf{RX}^{-1})(d\mathbf{X}) \hspace{-2cm} \nonumber\\
    &=&  \frac{\Gamma_{m}^{\beta}[a, -\kappa] \Gamma_{m}^{\beta}[(m-1)\beta/2+1]}
  {\Gamma_{m}^{\beta}[a + (m-1)\beta/2+1, -\kappa]}C_{\kappa}^{\beta}(\mathbf{R}),
\end{eqnarray}
for $\re(a) > (m-1)\beta/2 + k_{1}$. And
\begin{eqnarray}
  \int_{\mathbf{0} < \mathbf{X} < \mathbf{I}} |\mathbf{X}|^{a-(m-1)\beta/2-1}
    C_{\kappa}^{\beta}(\mathbf{XR})(d\mathbf{X}) \hspace{-2cm} \nonumber\\
    &=& \frac{\Gamma_{m}^{\beta}[a, \kappa]\Gamma_{m}^{\beta}[(m-1)\beta/2+1]}
    {\Gamma_{m}^{\beta}[a + (m-1)\beta/2+1, \kappa]} C_{\kappa}^{\beta}(\mathbf{R}),
\end{eqnarray}
for $\re(a) > (m-1)\beta/2 - k_{m}$.
\end{corollary}

Similarly, taking $a=(m-1)\beta/2 + 1 > (m-1)\beta/2 - k_{m}$ in (\ref{jpeq3}), we
have the following result.

\begin{corollary}\label{coro7}
If $\mathbf{R} \in \mathfrak{S}_{m}^{\beta,\mathfrak{C}}$, then
\begin{eqnarray}
  \int_{\mathbf{0} < \mathbf{X} < \mathbf{I}} |\mathbf{I} - \mathbf{X}|^{b-(m-1)\beta/2-1}
    C_{\kappa}^{\beta}(\mathbf{XR})(d\mathbf{X}) \hspace{-2cm} \nonumber\\
    &=& \frac{\Gamma_{m}^{\beta}[(m-1)\beta/2 + 1, \kappa]\Gamma_{m}^{\beta}[b]}{\Gamma_{m}^{\beta}[(m-1)\beta/2 + 1 + b, \kappa]}
    C_{\kappa}^{\beta}(\mathbf{R}),
\end{eqnarray}
for $\re(b) > (m-1)\beta/2$.
\end{corollary}

\section{Hypergeometric functions}\label{sec4}

In this section, we study diverse integral properties of hypergeometric functions for
normed division algebras. First, let us consider the following definition.

Fix complex numbers $a_{1}, \dots, a_{p}$ and $b_{1}, \dots, b_{q}$, and for all $1
\leq i \leq q$ and $1 \leq j \leq m$ do not allow $-b_{i} + (j-1)\beta/2$ to be a
nonnegative integer. Then the \emph{hypergeometric function with one matrix argument}
${}_{p}F_{q}^{\beta}$ is defined to be the real-analytic function on
$\mathfrak{S}_{m}^{\beta}$ given by the series
\begin{equation}\label{fhoa}
    {}_{p}F_{q}^{\beta}(a_{1}, \dots, a_{p};b_{1}, \dots, b_{q}; \mathbf{X}) = \sum_{k=0}^{\infty}\sum_{\kappa}
  \frac{[a_{1}]_{\kappa}^{\beta} \cdots [a_{p}]_{\kappa}^{\beta}}{[b_{1}]_{\kappa}^{\beta} \cdots
  [b_{q}]_{\kappa}^{\beta}} \ \frac{C_{\kappa}^{\beta}(\mathbf{X})}{k!}.
\end{equation}

Some known properties are, see \citet[Section 6, pp. 803-810]{gr:87}:

\emph{Convergence hypergeometric functions}.
\begin{enumerate}
  \item If $p \leq q$ then the hypergeometric series (\ref{fhoa}) converges absolutely
    for all $\mathbf{X} \in \mathfrak{S}_{m}^{\beta}$.
  \item If $p = q+1$ then the series (\ref{fhoa}) converges absolutely for $||\mathbf{X}||=
    \max\{|\lambda_{i}|: i = 1, \dots, m \} < 1$, and diverges for $||\mathbf{X}|| >
    1$, where $\lambda_{1}, \dots \lambda_{m}$ are the $i$-th eigenvalues of $\mathbf{X}
    \in \mathfrak{S}_{m}^{\beta}$.
  \item If $p > q$ then the series (\ref{fhoa}) diverges unless it terminates.
\end{enumerate}
For all $\mathbf{X} \in \mathfrak{S}_{m}^{\beta}$; indeed, for all $\mathbf{X} \in
\mathfrak{S}_{m}^{\beta, \mathfrak{C}}$. This is characteristic of the general
situation when $p \leq q$.
\begin{equation}\label{fh00}
    {}_{0}F_{0}^{\beta}(\mathbf{X}) = \sum_{k = 0}^{\infty}\sum_{\kappa}
    \frac{C_{\kappa}^{\beta}(\mathbf{X})}{k!} =  \sum_{k = 0}^{\infty} \frac{(\tr
    \mathbf{X})^{k}}{k!} = \etr\{\mathbf{X}\},
\end{equation}
If $\re(a) > (m-1)\beta/2$, and $||\mathbf{X}||< 1$,
\begin{eqnarray}\label{fh10}
  {}_{1}F_{0}^{\beta}(a;\mathbf{X}) &=& \frac{1}{\Gamma_{m}^{\beta}[a]}\int_{\mathbf{Y} \in
        \mathfrak{P}_{m}^{\beta}} \etr\{-(\mathbf{I} - \mathbf{X})\mathbf{Y}\}
        |\mathbf{Y}|^{a-(m-1)\beta/2-1} (d\mathbf{Y}) \\
    &=& |\mathbf{I} - \mathbf{X}|^{-a} \nonumber
\end{eqnarray}
gives the full analytic continuation of ${}_{1}F_{0}^{\beta}(a; \cdot)$ to any
simply-connected domain in $\mathfrak{S}_{m}^{\beta,\mathfrak{C}}$. The right side is
determined by the principal branch of the argument. The fact that the hypergeometric
series  ${}_{1}F_{0}^{\beta}$ has $\{\mathbf{X} \in \mathfrak{S}_{m}^{\beta}:
||\mathbf{X}|| < 1\}$ as its domain of convergence is characteristic of
${}_{p+1}F_{p}^{\beta}$ for all $p \geq 0$.

Let $\re(c) > \re(a) + (m-1)\beta/2 >(m-1)\beta$ and $||\mathbf{X}||< 1$. Then
$$
  {}_{p+1}F_{q+1}^{\beta}(a_{1}, \dots,a_{p},a;b_{1}, \dots, b_{q},c; \mathbf{X})
    \hspace{5cm}
$$\vspace{-0.5cm}
\begin{eqnarray}
     \phantom{******}&=& \frac{1}{\mathcal{B}_{m}^{\beta}[a,c-a]}
    \int_{\mathbf{0}<\mathbf{Y}<\mathbf{I}_{m}}
    {}_{p}F_{q}^{\beta}(a_{1} \cdots a_{p};b_{1} \cdots b_{q}; \mathbf{XY})
    \nonumber\\\label{fhp1q1}
    & & \qquad \times |\mathbf{Y}|^{a-(m-1)\beta/2-1} |\mathbf{I}-\mathbf{Y}|^{c-a-(m-1)\beta/2-1}(d\mathbf{Y}),
\end{eqnarray}
for $p = q+1$. In particular, for $p = 2$, we have the Euler formula
$$
    {}_{2}F_{1}^{\beta}(a_{1},a;c; \mathbf{X})=
    \frac{1}{\mathcal{B}_{m}^{\beta}[a,c-a]}\hspace{6cm}
$$ \vspace{-0.5cm}
\begin{eqnarray}
    \phantom{**********} && \times \int_{\mathbf{0}<\mathbf{Y}<\mathbf{I}_{m}}
    {}_{1}F_{0}^{\beta}(a_{1};\mathbf{XY})|\mathbf{Y}|^{a-(m-1)\beta/2-1}
    \nonumber\\\label{fh21}
    && \phantom{**********}\times |\mathbf{I}-\mathbf{Y}|^{c-a-(m-1)\beta/2-1}(d\mathbf{Y}).
\end{eqnarray}
for arbitrary $a_{1}$, $\re(c) > \re(a) + (m-1)\beta/2 >(m-1)\beta$ and
$||\mathbf{X}||< 1$.

\begin{remark}
Observe that (\ref{fhp1q1}) (and of course (\ref{fh21}), too) is a consequence of
(\ref{jpeq3}). And then the condition over $a$ and $c$ must be $\re(c) > \re(a) +
(m-1)\beta/2 >(m-1)\beta - k_{m}$.
\end{remark}

\emph{Laplace transform of hypergeometric functions}. Assume $p \leq q$, $\re(a) >
(m-1)\beta/2$ and $\mathbf{U} \in \mathfrak{S}_{m}^{\beta}$. Then
\begin{eqnarray}
    \int_{\mathbf{X} \in \mathfrak{P}_{m}^{\beta}}
    \etr\{-\mathbf{XZ}\}{}_{p}F_{q}^{\beta}(a_{1} \cdots a_{p};b_{1} \cdots b_{q};
    \mathbf{XU})|\mathbf{X}|^{a-(m-1)\beta/2-1} (d\mathbf{X}) \hspace{-0.5cm}\nonumber\\ \label{tlhfoa}
    \hspace{-0.5cm}= |\mathbf{Z}|^{-a} \Gamma_{m}^{\beta}[a] \ {}_{p+1}F_{q}^{\beta}(a_{1}
    \cdots a_{p},a;b_{1} \cdots b_{q}; \mathbf{UZ}^{-1}).
\end{eqnarray}
When $p < q$, the integral in (\ref{tlhfoa}) converges absolutely for all $\mathbf{Z}
\in \mathbf{\Phi}$. When $p = q$, the integral converges absolutely for all
$\mathbf{Z} \in \mathfrak{S}_{m}^{\beta, \mathfrak{C}}$, such that $
||(\re(\mathbf{Z}))^{-1} ||< 1.$

\begin{remark}\label{remtlhf}
\citet[p. 485]{h:55} wrote...``The convergence of the integral (\ref{tlhfoa})
requires at least $\re(a) > (m-1)\beta/2$..."However, \citet{gr:87} based their proof
of (\ref{tlhfoa}) in their Theorem 5.9 which is valid for $\re(a) > (m-1)\beta/2 -
k_{m}$. The same remak should be considered in the versions of \citet{c:63},
\citet{m:82}, \citet{rva:05b} and \citet{lx:09}.
\end{remark}

Similarly, the \emph{hypergeometric function of two matrix arguments}
${}_{p}F_{q}^{(m),\beta}$ is defined to be the real-analytic function on
$\mathfrak{S}_{m}^{\beta}$ given by the series%
{\small
\begin{equation}\label{fhta}
    {}_{p}F_{q}^{(m),\beta}(a_{1}, \dots, a_{p};b_{1}, \dots, b_{q}; \mathbf{X}, \mathbf{Y})
    = \sum_{k=0}^{\infty}\sum_{\kappa} \frac{[a_{1}]_{\kappa}^{\beta} \cdots
    [a_{p}]_{\kappa}^{\beta}}{[b_{1}]_{\kappa}^{\beta} \cdots [b_{q}]_{\kappa}^{\beta}}
    \ \frac{C_{\kappa}^{\beta}(\mathbf{X})C_{\kappa}^{\beta}(\mathbf{Y})}{k!
    \ C_{\kappa}^{\beta}(\mathbf{I})}.
\end{equation}}
Some basic properties of (\ref{fhta}) are shown below, see \citet{gr:89}.

\emph{Convergence hypergeometric functions with two matrix arguments}.
\begin{enumerate}
    \item If $p \leq q$ then the hypergeometric series (\ref{fhta}) converges absolutely
        for all $\mathbf{X}$ and $\mathbf{Y} \in \mathfrak{S}_{m}^{\beta}$.
    \item If $p = q+1$ then the series (\ref{fhta}) converges absolutely for $||\mathbf{X}||\cdot ||\mathbf{Y}||
        < 1$, and diverges for $||\mathbf{X}||\cdot ||\mathbf{Y}|| > 1$.
\end{enumerate}
Also
\begin{eqnarray}
    \int_{\mathbf{H} \in \mathfrak{U}^{\beta}(m)} {}_{p}F_{q}^{\beta}(a_{1}, \dots,
    a_{p};b_{1}, \dots, b_{q};
    \mathbf{XHYH}^{*})(d\mathbf{H})\hspace{2.5cm}\\\label{fhta1}
    = {}_{p}F_{q}^{(m),\beta}(a_{1}, \dots, a_{p};b_{1}, \dots, b_{q}; \mathbf{X},
    \mathbf{Y}).
\end{eqnarray}
In particular
\begin{eqnarray}
    \int_{\mathbf{H} \in \mathfrak{U}^{\beta}(m)} {}_{0}F_{0}^{\beta}(
    \mathbf{XHYH}^{*})(d\mathbf{H}) &=&  \int_{\mathbf{H} \in \mathfrak{U}^{\beta}(m)}
    \etr\{\mathbf{XHYH}^{*}\}(d\mathbf{H})\\\label{exp}
    &=& {}_{0}F_{0}^{(m),\beta}(\mathbf{X}, \mathbf{Y}).
\end{eqnarray}

\emph{Laplace transform of hypergeometric functions with two matrix arguments}.

Assume $p \leq q$, $\re(a) > (m-1)\beta/2$ and $\mathbf{U} \in
\mathfrak{S}_{m}^{\beta}$. Then
$$
  \int_{\mathbf{X} \in \mathfrak{P}_{m}^{\beta}}
  \etr\{-\mathbf{XZ}\}{}_{p}F_{q}^{(m),\beta}(a_{1} \cdots a_{p};b_{1} \cdots b_{q};
  \mathbf{XU}, \mathbf{Y})|\mathbf{X}|^{a-(m-1)\beta/2-1} (d\mathbf{X})
$$\vspace{-0.5cm}
\begin{equation} \label{tlhfta}
    \phantom{****} = |\mathbf{Z}|^{-a}\Gamma_{m}^{\beta}[a] \ {}_{p+1}F_{q}^{(m),\beta}
    (a_{1} \cdots a_{p},a;b_{1} \cdots b_{q}; \mathbf{UZ}^{-1},\mathbf{Y}).
\end{equation}
When $p < q$, the integral in (\ref{tlhfta}) converges absolutely for all $\mathbf{Z}
\in \mathbf{\Phi}$ and $\mathbf{Y} \in \mathfrak{S}_{m}^{\beta}$. When $p = q$, the
integral converges absolutely for all $\mathbf{Z}$ and $\mathbf{Y} \in
\mathfrak{S}_{m}^{\beta, \mathfrak{C}}$, such that $ ||(\re(\mathbf{Z}))^{-1} ||\cdot
||\mathbf{Y}||< 1.$ Remark similar tov\ref{remtlhf} should be considered for
(\ref{tlhfta}).

We now propose further integral properties of hypergeometric functions for normed
division algebras. The first result is the inverse Laplace transformation. In the
real case, this result was obtained by \citet{h:55}, \citet{c:63}, \citet{j:64} and
\citet[p. 261]{m:82}, and in the complex case by \citet[p. 370]{m:97}. Let us first
consider the following extension of similar results discussed in \citet{c:63}, see
also \citet[p. 253]{m:82}.

\begin{lemma}\label{lemitl}
Assume that $\mathbf{Z}$ and  $\mathbf{X} \in \mathfrak{S}_{m}^{\beta,
\mathfrak{C}}$, $\mathbf{U} \in \mathfrak{S}_{m}^{\beta}$ and $\re(a) > a_{0}$. Then
$$
   \frac{2^{m(m-1)\beta/2}}{(2\pi i)^{m(m-1)\beta/2+ m}} \int_{\mathbf{Z}
   - \mathbf{Z}_{0} \in \mathbf{\Phi}} \etr\{\mathbf{XZ}\} |\mathbf{Z}|^{-a}
    C_{\kappa}^{\beta}(\mathbf{UZ}^{-1})(d\mathbf{Z})\hspace{3cm}
$$\vspace{-0.5cm}
\begin{equation}\label{eqlemitl}\hspace{4cm}
   = \frac{1}{\Gamma_{m}^{\beta}[a,\kappa]}|\mathbf{X}|^{a-(m-1)\beta/2-1}C_{\kappa}^{\beta}(\mathbf{XU}),
\end{equation}
where $\mathbf{Z}_{0} \in \mathfrak{P}_{m}^{\beta}$.
\end{lemma}

\begin{theorem}\label{tilt} Assume that $\mathbf{Z}$ and $\mathbf{X} \in
\mathfrak{S}_{m}^{\beta, \mathfrak{C}}$, $\mathbf{U} \in \mathfrak{S}_{m}^{\beta}$
and $\re(b) > b_{0}$. Then
$$
  \frac{\Gamma_{m}^{\beta}[b] 2^{m(m-1)\beta/2}}{(2\pi i)^{m(m-1)\beta/2+ m}} \int_{\mathbf{Z}
  - \mathbf{Z}_{0} \in \mathbf{\Phi}} \etr\{\mathbf{XZ}\} |\mathbf{Z}|^{-b}
  {}_{p}F_{q}^{\beta}(a_{1},\dots,a_{p}; b_{1}, \dots,b_{q}; \mathbf{UZ}^{-1})(d\mathbf{Z})
$$
\begin{equation}\label{eqtil1}
    \hspace{3cm}= |\mathbf{X}|^{b-(m-1)\beta/2-1}{}_{p}F_{q+1}^{\beta}(a_{1},\dots,a_{p};
    b_{1}, \dots,b_{q},b; \mathbf{XU}),
\end{equation}
and if $\mathbf{Y} \in \mathfrak{S}_{m}^{\beta}$
\begin{eqnarray*}\frac{\Gamma_{m}^{\beta}[b] 2^{m(m-1)\beta/2}}{(2\pi i)^{m(m-1)\beta/2+ m}} \int_{\mathbf{Z}
  -  \mathbf{Z}_{0} \in \mathbf{\Phi}} \etr\{\mathbf{XZ}\} |\mathbf{Z}|^{-b} \hspace{4cm}\\
  \times \ {}_{p}F_{q}^{(m),\beta}(a_{1},\dots,a_{p}; b_{1}, \dots,b_{q}; \mathbf{UZ}^{-1},\mathbf{Y})(d\mathbf{Z})
\end{eqnarray*}
\begin{equation}\label{eqtil2}
    \hspace{2.5cm}= |\mathbf{X}|^{b-(m-1)\beta/2-1}{}_{p}F_{q+1}^{(m),\beta}(a_{1},\dots,a_{p};
    b_{1}, \dots,b_{q},b; \mathbf{XU},\mathbf{Y}),
\end{equation}
where $\mathbf{Z}_{0} \in \mathfrak{P}_{m}^{\beta}$.
\end{theorem}
\begin{proof}
Proof  of both (\ref{eqtil1}) and (\ref{eqtil2}) follows by expanding the
${}_{p}F_{q}^{\beta}$ and ${}_{p}F_{q}^{(m),\beta}$ functions in the integrands and
integrating term by term using (\ref{eqlemitl}).
\end{proof}

\begin{theorem}\label{teofh11}
The ${}_{1}F_{1}^{\beta}$ function has the integral representation
$$
  {}_{1}F_{1}^{\beta}(a,c;\mathbf{X}) = \frac{1}{\mathcal{B}_{m}^{\beta}[a,c-a]}
  \int_{\mathbf{0} < \mathbf{Y} < \mathbf{I}} \etr\{\mathbf{XY}\}
  |\mathbf{Y}|^{a-(m-1)\beta/2-1} \qquad\qquad\qquad
$$\vspace{-0.5cm}
\begin{equation}\label{fh11}
    \qquad\qquad\times |\mathbf{I} - \mathbf{Y}|^{c-a-(m-1)\beta/2-1}(d\mathbf{Y}),
\end{equation}
valid for $\re(c) > \re(a) + (m-1)\beta/2 >(m-1)\beta - k_{m}$ and all $\mathbf{X}
\in \mathfrak{S}_{m}^{\beta,\mathfrak{C}}$.
\end{theorem}
\begin{proof}
The desired result is obtained by expanding $\etr\{\mathbf{XY}\}$ using (\ref{jpexp})
and integrating term by term using (\ref{jpeq3}).
\end{proof}

The generalised Kummer and Euler relations are given in the following result.

\begin{theorem}\label{teoKEr}
\begin{equation}\label{kr}
    {}_{1}F_{1}^{\beta}(a,c;\mathbf{X}) =
    \etr\{\mathbf{X}\}{}_{1}F_{1}^{\beta}(c-a,c;-\mathbf{X}),
\end{equation}
for $\mathbf{X} \in \mathfrak{S}_{m}^{\beta,\mathfrak{C}}$. And
\begin{eqnarray}\label{er1}
  {}_{2}F_{1}^{\beta}(a,b;c;\mathbf{X}) &=& |\mathbf{I} - \mathbf{X}|^{b}
  {}_{2}F_{1}^{\beta}(c-a,b;c;-\mathbf{X}(\mathbf{I} - \mathbf{X})^{-1})
  \\\label{er2}
   &=& |\mathbf{I} - \mathbf{X}|^{c-a-b}
  {}_{2}F_{1}^{\beta}(c-a,c-b;c;\mathbf{X})
\end{eqnarray}
for $||\mathbf{X}||< 1$.
\end{theorem}

\begin{remark}
Observe that, for ${}_{1}F_{0}^{\beta}(a;\mathbf{X})$ the condition $\re(a) >
(m-1)\beta/2$ over $a$ is determined by its integral representation (\ref{fh10}).
However ${}_{1}F_{0}^{\beta}(a;\mathbf{X})$ is easily seen to be analytic for all $a$
and $||\mathbf{X}||< 1$, see \citet[p. 486]{h:55}. Similarly, the conditions $\re(c)
> \re(a) + (m-1)\beta/2 >(m-1)\beta - k_{m}$ over $a$ and $c$ given in Theorem
\ref{teofh11}, valid for Theorem \ref{teoKEr} (\ref{kr}) too, are determined by the
existence of the integral (\ref{fh11}). However, these conditions can be extended to
other possible values if we use the inverse Laplace transformation to define
${}_{1}F_{1}^{\beta}(a,c;\mathbf{X})$. In this case
${}_{1}F_{1}^{\beta}(a,c;\mathbf{X})$ is valid for the arbitrary complex $a$, $\re(c)
> (m - 1)\beta/2$ and $\mathbf{X} \in \mathfrak{S}_{m}^{\beta,\mathfrak{C}}$, see
\citet[p. 487]{h:55}. Also, the conditions $\re(c) > \re(a) + (m-1)\beta/2
>(m-1)\beta - k_{m}$ over $a$ and $c$ for (\ref{fh21}) and Theorem
\ref{teoKEr}(\ref{er1}) and (\ref{er2}) are determined by the absolutely convergence
of the integral (\ref{fh21}). Again, these conditions about $a$ and $c$ can be
extended to other possible values using the inverse Laplace transformation and the
results for ${}_{1}F_{1}^{\beta}(a,c;\mathbf{X})$ obtained as described before, see
\citet[p. 489]{h:55}. Finally, let us take into account that, for any analysis if the
integral representation of ${}_{1}F_{1}^{\beta}(a,c;\mathbf{X})$ or
${}_{2}F_{1}^{\beta}(a,b;c;\mathbf{X})$  is not used explicitly, then the extended
conditions for $a$ and $c$ could be considered.

\end{remark}

\begin{theorem}\label{teofh01}
Let $\mathbf{X} \in {\mathcal L}^{\beta}_{m,n}$ and $\mathbf{H} =
(\mathbf{H}_{1}|\mathbf{H}_{2}) \in \mathfrak{U}^{\beta}(m)$, $\mathbf{H}_{1} \in
\mathcal{V}_{m,n}^{\beta}$. Then
\begin{eqnarray}\label{eqfh011}
    {}_{0}F_{1}^{\beta} (\beta n/2; \beta^{2} \mathbf{XX}^{*}/4) &=&
    \int_{\mathbf{H} \in \mathfrak{U}^{\beta}(m)} \etr(\beta
    \mathbf{XH}_{1})(d\mathbf{H})\\ \label{eqfh012}
    &=& \int_{\mathbf{H}_{1} \in \mathcal{V}_{m,n}^{\beta}} \etr(\beta
    \mathbf{XH}_{1})(d\mathbf{H}_{1})
\end{eqnarray}
\end{theorem}
\begin{proof}
The proof is analogous to that given in the real case for \citet[Theorem 7.4.1]{m:82}
and in the quaternion case by \citet{lx:09}. Alternative proofs can be established in
an analogous form to those given by \citet{j:61} and \citet[p. 494-495]{h:55}. For
(\ref{eqfh012}) it might be necessary to consider Lemma 9.5.3, p. 397 in
\citet{m:82}.
\end{proof}

On the basis of Theorem \ref{teo1}, we now discuss diverse integral properties of
generalised hypergeometric functions, which contain as particular cases many of the
results established above.

\begin{theorem}\label{teofhg1}
Assume $p \leq q$ and $\re(a) > (m-1)\beta/2 - k_{m}$ and $\mathbf{U} \in
\mathfrak{S}_{m}^{\beta}$. Then for $\vartheta = \int_{0}^{\infty}f(z) z^{am+k-1} dz
< \infty$,
\begin{eqnarray}
  \int_{\mathbf{X} \in \mathfrak{P}_{m}^{\beta}}
  f(\tr \mathbf{XZ}){}_{p}F_{q}^{\beta}(a_{1} \cdots a_{p};b_{1} \cdots b_{q};
  \mathbf{XU})|\mathbf{X}|^{a-(m-1)\beta/2-1} (d\mathbf{X})\nonumber\\ \label{tlghfoa}
  = |\mathbf{Z}|^{-a} \Gamma_{m}^{\beta}[a] \ \sum_{k=0}^{\infty}\sum_{\kappa}
  \frac{[a_{1}]_{\kappa}^{\beta} \cdots [a_{p}]_{\kappa}^{\beta}[a]_{\kappa}^{\beta}}
  {[b_{1}]_{\kappa}^{\beta} \cdots
  [b_{q}]_{\kappa}^{\beta}} \ \frac{C_{\kappa}^{\beta}(\mathbf{UZ}^{-1})}{\Gamma[am+k]k!} \cdot \vartheta.
\end{eqnarray}
When $p < q$, the integral in (\ref{tlghfoa}) converges absolutely for all
$\mathbf{Z} \in \mathbf{\Phi}$. When $p = q$, the integral converges absolutely for
all $\mathbf{Z} \in \mathfrak{S}_{m}^{\beta, \mathfrak{C}}$, such that $
||(\re(\mathbf{Z}))^{-1} ||< 1.$\\%
Similarly, let $p \leq q$, $\re(a) > (m-1)\beta/2 - k_{m}$ and $\mathbf{U} \in
\mathfrak{S}_{m}^{\beta}$. Then
$$
  \int_{\mathbf{X} \in \mathfrak{P}_{m}^{\beta}}
  f(\tr \mathbf{XZ}){}_{p}F_{q}^{(m),\beta}(a_{1} \cdots a_{p};b_{1} \cdots b_{q};
  \mathbf{XU}, \mathbf{Y})|\mathbf{X}|^{a-(m-1)\beta/2-1} (d\mathbf{X})
  $$\vspace{-0.5cm}
\begin{equation} \label{tlghfta}
    = |\mathbf{Z}|^{-a}\Gamma_{m}^{\beta}[a] \sum_{k=0}^{\infty}\sum_{\kappa}
    \frac{[a_{1}]_{\kappa}^{\beta} \cdots
    [a_{p}]_{\kappa}^{\beta}[a]_{\kappa}^{\beta}}{[b_{1}]_{\kappa}^{\beta} \cdots [b_{q}]_{\kappa}^{\beta}}
    \ \frac{C_{\kappa}^{\beta}(\mathbf{UZ}^{-1})C_{\kappa}^{\beta}(\mathbf{Y})}{\Gamma[am+k]k!
    \ C_{\kappa}^{\beta}(\mathbf{I})} \cdot \vartheta.
\end{equation}
When $p < q$, the integral in (\ref{tlghfta}) converges absolutely for all
$\mathbf{Z} \in \mathbf{\Phi}$ and $\mathbf{Y} \in \mathfrak{S}_{m}^{\beta}$. When $p
= q$, the integral converges absolutely for all $\mathbf{Z}$ and $\mathbf{Y} \in
\mathfrak{S}_{m}^{\beta, \mathfrak{C}}$, such that $ ||(\re(\mathbf{Z}))^{-1} ||\cdot
||\mathbf{Y}||< 1.$
\end{theorem}

Observe that if $f(\tr \mathbf{XZ}) = \etr\{- \mathbf{XZ}\}$ in Theorem \ref{teofhg1}
then we obtain (\ref{tlhfoa}) and (\ref{tlhfta}).

\begin{theorem}\label{teofhg2}
Assume $p \leq q$ and $\re(a) > (m-1)\beta/2 - k_{m}$ and $\mathbf{U} \in
\mathfrak{S}_{m}^{\beta}$. Then for $\gamma = \int_{0}^{\infty}f(z) z^{am+k-1} dz <
\infty$,
$$
  \int_{\mathbf{X} \in \mathfrak{P}_{m}^{\beta}}
  f(\tr \mathbf{XZ}){}_{p}F_{q}^{\beta}(a_{1} \cdots a_{p};b_{1} \cdots b_{q};
  \mathbf{X}^{-1}\mathbf{U})|\mathbf{X}|^{a-(m-1)\beta/2-1} (d\mathbf{X})\hspace{2cm}
$$\vspace{-0.5cm}
\begin{eqnarray}
    &=& |\mathbf{Z}|^{-a} \Gamma_{m}^{\beta}[a] \ \sum_{k=0}^{\infty}\sum_{\kappa}
    \frac{(-1)^{k}[a_{1}]_{\kappa}^{\beta} \cdots [a_{p}]_{\kappa}^{\beta}[a]_{\kappa}^{\beta}}
    {[b_{1}]_{\kappa}^{\beta} \cdots [b_{q}]_{\kappa}^{\beta}[-a +(m-1)\beta/2 +1]_{\kappa}^{\beta}}
    \nonumber\\ \label{tlghfoainv}
    & & \hspace{8cm}\times \ \frac{C_{\kappa}^{\beta}(\mathbf{UZ})}{\Gamma[am+k]k!} \cdot \gamma.
\end{eqnarray}
When $p < q$, the integral in (\ref{tlghfoainv}) converges absolutely for all
$\mathbf{Z} \in \mathbf{\Phi}$. When $p = q$, the integral converges absolutely for
all $\mathbf{Z} \in \mathfrak{S}_{m}^{\beta, \mathfrak{C}}$, such that $
||(\re(\mathbf{Z}))^{-1} ||< 1.$\\%
Similarly, let $p \leq q$, $\re(a) > (m-1)\beta/2 - k_{m}$ and $\mathbf{U} \in
\mathfrak{S}_{m}^{\beta}$. Then
$$
  \int_{\mathbf{X} \in \mathfrak{P}_{m}^{\beta}}
  f(\tr \mathbf{XZ}){}_{p}F_{q}^{(m),\beta}(a_{1} \cdots a_{p};b_{1} \cdots b_{q};
  \mathbf{X}^{-1}\mathbf{U}, \mathbf{Y})|\mathbf{X}|^{a-(m-1)\beta/2-1} (d\mathbf{X})
  $$\vspace{-0.5cm}
\begin{eqnarray}
    &=& |\mathbf{Z}|^{-a}\Gamma_{m}^{\beta}[a] \sum_{k=0}^{\infty}\sum_{\kappa}
    \frac{(-1)^{k}[a_{1}]_{\kappa}^{\beta} \cdots
    [a_{p}]_{\kappa}^{\beta}[a]_{\kappa}^{\beta}}{[b_{1}]_{\kappa}^{\beta} \cdots
    [b_{q}]_{\kappa}^{\beta} [-a +(m-1)\beta/2 +1]_{\kappa}^{\beta}}\nonumber\\ \label{tlghftainv}
    && \hspace{7cm}\times \ \frac{C_{\kappa}^{\beta}(\mathbf{UZ})C_{\kappa}^{\beta}(\mathbf{Y})}{\Gamma[am+k]k!
    \ C_{\kappa}^{\beta}(\mathbf{I})} \cdot \gamma.
\end{eqnarray}
When $p < q$, the integral in (\ref{tlghftainv}) converges absolutely for all
$\mathbf{Z} \in \mathbf{\Phi}$ and $\mathbf{Y} \in \mathfrak{S}_{m}^{\beta}$. When $p
= q$, the integral converges absolutely for all $\mathbf{Z}$ and $\mathbf{Y} \in
\mathfrak{S}_{m}^{\beta, \mathfrak{C}}$, such that $ ||(\re(\mathbf{Z}))^{-1} ||\cdot
||\mathbf{Y}||< 1.$
\end{theorem}

Now, let us define $f(\tr \mathbf{XZ}) = \etr\{- \mathbf{XZ}\}$ in Theorem
\ref{teofhg2}, then we obtain:
\begin{eqnarray}
    \int_{\mathbf{X} \in \mathfrak{P}_{m}^{\beta}}
    \etr\{-\mathbf{XZ}\}{}_{p}F_{q}^{\beta}(a_{1} \cdots a_{p};b_{1} \cdots b_{q};
    \mathbf{X}^{-1}\mathbf{U})|\mathbf{X}|^{a-(m-1)\beta/2-1} (d\mathbf{X})\hspace{0.5cm}\nonumber\\ \label{tlhfoainv}
    = |\mathbf{Z}|^{-a} \Gamma_{m}^{\beta}[a] \ {}_{p}F_{q+1}^{\beta}(a_{1} \cdots a_{p};b_{1} \cdots b_{q},-a +(m-1)\beta/2 +1;
    -\mathbf{UZ}).
\end{eqnarray}
and
$$
  \int_{\mathbf{X} \in \mathfrak{P}_{m}^{\beta}}
  \etr\{-\mathbf{XZ}\}{}_{p}F_{q}^{(m),\beta}(a_{1} \cdots a_{p};b_{1} \cdots b_{q};
  \mathbf{X}^{-1}\mathbf{U}, \mathbf{Y})|\mathbf{X}|^{a-(m-1)\beta/2-1} (d\mathbf{X})
  $$\vspace{-0.5cm}
\begin{equation} \label{tlhftainv}
    = |\mathbf{Z}|^{-a}\Gamma_{m}^{\beta}[a] \ {}_{p}F_{q+1}^{(m),\beta}
    (a_{1} \cdots a_{p};b_{1} \cdots b_{q},-a +(m-1)\beta/2 +1;
    -\mathbf{UZ},\mathbf{Y}),
\end{equation}
where for both $\re(a) > (m-1)\beta/2 + k_{1}$.

Similar results to (\ref{tlhfoa}) and (\ref{tlhfta}) or (\ref{tlhfoainv}) and
(\ref{tlhftainv}) can be obtained from Corollaries \ref{coro3} and \ref{coro4}.

Now, we propose the incomplete gamma and beta functions for normed division algebras.

\begin{theorem}\label{teoigb1}
Let $\mathbf{\Lambda} \in \mathfrak{S}_{m}^{\beta,\mathfrak{C}}$ and $\mathbf{\Omega}
\in \mathbf{\Phi}$. Then
\begin{eqnarray}
   \int_{\mathbf{0}< \mathbf{X} < \mathbf{\Omega}} \etr\{-\mathbf{\Lambda}\mathbf{X}\}
   |\mathbf{X}|^{a-(m-1)\beta/2-1}(d\mathbf{X})\hspace{5cm} \nonumber\\\label{igf}
   =  \mathcal{B}_{m}^{\beta}[a, (m-1)\beta/2+1]
   |\mathbf{\Omega}|^{a} {}_{1}F_{1}^{\beta}(a; a + (m -1)\beta/2 +1; -\mathbf{\Omega\Lambda}),
\end{eqnarray}
for $\re(a) > (m-1)\beta/2 - k_{m}$. And, let $\mathbf{0} < \mathbf{\Xi} <
\mathbf{I}$, then
$$
   \int_{\mathbf{0} < \mathbf{Y} < \mathbf{\Xi}} |\mathbf{Y}|^{a-(m-1)\beta/2-1}
   |\mathbf{I} - \mathbf{Y}|^{b-(m-1)\beta/2-1}(d\mathbf{Y}) = \mathcal{B}_{m}^{\beta}
   [a, (m-1)\beta/2+1]
$$ \vspace{-0.5cm}
\begin{equation}\label{ibf}
   \times \ |\mathbf{\Xi}|^{a} {}_{2}F_{1}^{\beta}(a, -b + (m - 1)\beta/2 + 1;
   a + (m - 1)\beta/2 + 1; \mathbf{\Xi}),
\end{equation}
for $\re(a) > (m-1)\beta/2 - k_{m}$ and $\re(b) > (m-1)\beta/2$.
\end{theorem}
\begin{proof}
For (\ref{igf}), let us make the transformation $\mathbf{X} = \mathbf{\Omega}^{1/2}
\mathbf{R\Omega}^{1/2}$ and by applying Lemma \ref{lemhlt} we have, $(d\mathbf{X}) =
|\mathbf{\Omega}|^{(m+1)\beta/2 +1}(d\mathbf{R})$, with $\mathbf{0} < \mathbf{R} <
\mathbf{I}$. Then, expanding $\etr\{-\mathbf{\Lambda X}\}$ as a series of Jack
polynomials and integrating term by term using Corollary \ref{coro6}, the desired
result is obtained. Similarly, (\ref{ibf}) is proved by making the transformation
$\mathbf{Y} = \mathbf{\Xi}^{1/2}\mathbf{R \Xi}^{1/2}$ from where, applying the Lemma
\ref{lemhlt} we obtain that $(d\mathbf{X}) = |\mathbf{\Xi}|^{(m+1)\beta/2
+1}(d\mathbf{R})$, with $\mathbf{0} < \mathbf{R} < \mathbf{I}$, expanding
$|\mathbf{I} - \mathbf{X\Xi}|^{b-(m-1)\beta/2-1} = {}_{1}F_{0}^{\beta}(-b +
(m-1)\beta/2 + 1; \mathbf{X\Xi})$ and integrating term by term using Corollary
\ref{coro6}.
\end{proof}

\begin{theorem}\label{teoigb2}
Let $\mathbf{\Lambda} \in \mathfrak{S}_{m}^{\beta,\mathfrak{C}}$ and $\mathbf{\Omega}
\in \mathbf{\Phi}$. If $r = a-(m-1)\beta/2-1$ is a positive integer, then
\begin{eqnarray}
   \int_{\mathbf{X} > \mathbf{\Omega}} \etr\{-\mathbf{\Lambda}\mathbf{X}\}
   |\mathbf{X}|^{a-(m-1)\beta/2-1}(d\mathbf{X})\hspace{4cm} \nonumber\\\label{igf2}
   =  \Gamma_{m}^{\beta}[a]|\mathbf{\Lambda}|^{-a} \etr\{-\mathbf{\Lambda \Omega}\}
   \sum_{k =0}^{mr}\sum_{\kappa}{}^{*}
   \frac{C_{\kappa}^{\beta}(\mathbf{\Omega\Lambda})}{k!},
\end{eqnarray}
for $\re(a) > (m-1)\beta/2 + k_{1}$ and $\sum_{\kappa}^{*}$ denotes summation over
those partitions $\kappa = (k_{1}, \dots, k_{m})$ of $k$ with $k_{1} \leq r$.
\end{theorem}
\begin{proof}
Consider the transformation $\mathbf{X} = \mathbf{\Omega}^{1/2}(\mathbf{I} +
\mathbf{R}) \mathbf{\Omega}^{1/2}$ and applying Lemma \ref{lemhlt} we have,
$(d\mathbf{X}) = |\mathbf{\Omega}|^{(m+1)\beta/2 +1}(d\mathbf{R})$, with $\mathbf{R}
> \mathbf{0}$. Noting that $|\mathbf{I} + \mathbf{R}| = |\mathbf{R}||\mathbf{I} +
\mathbf{R}^{-1}|$ and expanding $|\mathbf{I} + \mathbf{R}^{-1}|^{a-(m-1)\beta/2-1}$
in terms of Jack polynomials, assuming that $r = a-(m-1)\beta/2-1$ is a positive
integer we obtain
\begin{eqnarray*}
  |\mathbf{I} + \mathbf{R}^{-1}|^{a-(m-1)\beta/2-1} &=& {}_{1}F_{0}^{\beta}(-a+(m-1)\beta/2+1; - \mathbf{R}^{-1}) \\
  &=& \sum_{k =0}^{mr}\sum_{\kappa}{}^{*}
   \frac{[-a+(m-1)\beta/2+1]_{\kappa}^{\beta} (-1)^{k}C_{\kappa}^{\beta}(\mathbf{\Omega\Lambda})}{k!}
\end{eqnarray*}
because $[-a+(m-1)\beta/2+1]_{\kappa}^{\beta} \equiv 0$ is any part of $\kappa$ that
is greater than $r$. The desired result is obtained by integrating term by term using
Corollary \ref{coro2}.
\end{proof}

We end this section with a some general results, which are useful in a variety of
situations, which enable us to transform the density function of a matrix $\mathbf{X}
\in \mathfrak{P}_{m}^{\beta}$ to the density function of its eigenvalues.

\begin{theorem}
Let $\mathbf{X} \in \mathfrak{P}_{m}^{\beta}$ be a random matrix with density
function $f(\mathbf{X})$. Then the joint density function of the eigenvalues
$\lambda_{1}, \dots, \lambda_{m}$ of $\mathbf{X}$ is
\begin{equation}\label{dfeig}
    \frac{\pi^{m^{2}\beta/2+ \varrho}}{\Gamma_{m}^{\beta}[m\beta/2]} \prod_{i <
    j}^{m}(\lambda_{i} -\lambda_{j})^{\beta}\int_{\mathbf{H} \in
    \mathfrak{U}^{\beta}(m)}f(\mathbf{HLH}^{*})(d\mathbf{H})
\end{equation}
where $\mathbf{L} = \diag(\lambda_{1}, \dots, \lambda_{m})$, $\lambda_{1}> \cdots >
\lambda_{m} > 0$, $\varrho$ is defined in Lemma \ref{lemsd} and $(d\mathbf{H})$ is
the normalised Haar measure.
\end{theorem}
\begin{proof}
The proof follows immediately from Lemma \ref{lemsd}.
\end{proof}

\section{Invariant polynomials}\label{sec5}

In this section, we extend many of the properties of a class of homogeneous
polynomials for normed division algebras of degrees $k$ and $t$ in the elements of
matrices $\mathbf{X}$ and $\mathbf{Y} \in \mathfrak{S}_{m}^{\beta}$, respectively,
see \citet{da:79}, \citet{da:80}, \citet{ch:80} and \citet{chd:86}; these are denoted
as $C_{\phi}^{[\beta]\kappa,\tau}(\mathbf{X},\mathbf{Y})$. These homogeneous
polynomials are invariant under the simultaneous transformations
$$
  \mathbf{X} \rightarrow \mathbf{U}^{*}\mathbf{XU}, \qquad \mathbf{Y} \rightarrow
  \mathbf{U}^{*}\mathbf{YU}, \qquad \mathbf{H} \in \mathfrak{U}^{\beta}(m).
$$
The most important relationship of these polynomials is
\begin{eqnarray}
    \int_{\mathbf{H} \in \mathfrak{U}^{\beta}(m)} C_{\kappa}^{\beta}(\mathbf{AH}^{*}\mathbf{XH})
    C_{\tau}^{\beta}(\mathbf{BH}^{*}\mathbf{YH})(d\mathbf{H}) \hspace{4cm}\nonumber\\\label{pi1}
    = \sum_{\phi \in \kappa.\tau}
    \frac{C_{\phi}^{[\beta]\kappa, \tau}(\mathbf{A}, \mathbf{B})
    C_{\phi}^{[\beta]\kappa, \tau}(\mathbf{X},
    \mathbf{Y})}{C_{\phi}^{\beta}(\mathbf{I})},
\end{eqnarray}
where $(d\mathbf{H})$ is the normalised Haar measure and $C_{\kappa}^{\beta}$,
$C_{\tau}^{\beta}$ and $C_{\phi}^{\beta}$ are Jack polynomials indexed by ordered
partitions $\kappa$, $\tau$ and $\phi$ of nonnegative integers $k$, $t$ and $f = k +
t$, respectively, into not more than $m$ parts. $\phi \in \kappa.\tau$ denotes the
irreducible representation of $GL(m,\mathfrak{F})$ indexed by $2\phi$ that occurs in
the decomposition of the Kronecker product $2\kappa \otimes 2\tau$ of the irreducible
representations indexed by $2\kappa$ and $2\tau$, see \citet{da:79} and
\citet{da:80}.

In a similar way to the case of Jack polynomials, let $\mathbf{A} =
\mathbf{A}^{*}\mathbf{A}$ and $\mathbf{B} = \mathbf{B}^{*}\mathbf{B}$. We opt for
convenience of notation rather than strict adherence to rigor, and write
$C_{\kappa}^{[\beta],\kappa,\tau}(\mathbf{XA}, \mathbf{YB})$ or
$C_{\kappa}^{[\beta],\kappa,\tau}(\mathbf{AX}, \mathbf{BY})$ rather than
$C_{\kappa}^{[\beta],\kappa,\tau}(\mathbf{AXA}^{*}, \mathbf{BYB}^{*})$, even though
$\mathbf{XA}$, $\mathbf{YB}$, $\mathbf{AX}$, or $\mathbf{BY}$ need not lie in
$\mathfrak{S}_{m}^{\beta}$.

Some of the elementary properties and results on invariant polynomials are extended
below:

\bigskip%
\textbf{Elementary properties of $\ C_{\phi}^{[\beta]\kappa, \tau}$}.%
\bigskip

Let $\mathbf{X}$ and $\mathbf{Y} \in \mathfrak{S}_{m}^{\beta}$, then
\begin{equation}\label{pi2}
    C_{\phi}^{[\beta]\kappa, \tau}(\mathbf{X}, \mathbf{X}) =
    \theta_{\phi}^{[\beta]\kappa, \tau} C_{\phi}^{\beta}(\mathbf{X}), \quad \mbox{where}
    \quad \theta_{\phi}^{[\beta]\kappa, \tau} = \frac{{}^{\beta}
    C_{\phi}^{\kappa, \tau}(\mathbf{I}, \mathbf{I})}{C_{\phi}^{\beta}(\mathbf{I})}.
\end{equation}
\begin{equation}\label{pi3}
    C_{\phi}^{[\beta]\kappa, \tau}(\mathbf{X}, \mathbf{Y}) =
    \left\{
      \begin{array}{ll}
        \displaystyle\frac{\theta_{\phi}^{[\beta]\kappa, \tau} C_{\phi}^{\beta}(\mathbf{I})}{C_{\kappa}^{\beta}(\mathbf{I})}
        C_{\kappa}^{\beta}(\mathbf{X}), & \hbox{ for } \mathbf{Y} = \mathbf{I};\\
        \displaystyle\frac{\theta_{\phi}^{[\beta]\kappa, \tau} C_{\phi}^{\beta}(\mathbf{I})}{C_{\tau}^{\beta}(\mathbf{I})}
        C_{\tau}^{\beta}(\mathbf{Y}), & \hbox{ for } \mathbf{X} = \mathbf{I}.
      \end{array}
    \right.
\end{equation}
\begin{equation}\label{pi4}
    C_{\phi}^{[\beta]\kappa, 0}(\mathbf{X}, \mathbf{Y}) =
     C_{\kappa}^{\beta}(\mathbf{X}), \quad \mbox{and}
    \quad C_{\phi}^{[\beta]0, \tau}(\mathbf{X}, \mathbf{Y}) =
     C_{\tau}^{\beta}(\mathbf{Y}).
\end{equation}
\begin{equation}\label{pi5}
    C_{\kappa}^{\beta}(\mathbf{X}) C_{\tau}^{\beta}(\mathbf{Y}) = \sum_{\phi \in \kappa.\tau}
    \theta_{\phi}^{[\beta]\kappa, \tau} C_{\phi}^{[\beta]\kappa, \tau}(\mathbf{X},
    \mathbf{Y}),
\end{equation}
therefore,
\begin{equation}\label{pi6}
    (\tr \mathbf{X})^{k} (\tr \mathbf{Y})^{t} = \sum_{\kappa, \tau;\phi \in \kappa.\tau}
    \theta_{\phi}^{[\beta]\kappa, \tau} C_{\phi}^{[\beta]\kappa, \tau}(\mathbf{X}, \mathbf{Y}).
\end{equation}
From (\ref{pi2}) and (\ref{pi5})
\begin{equation}\label{pi7}
    C_{\kappa}^{\beta}(\mathbf{X}) C_{\tau}^{\beta}(\mathbf{X}) = \sum_{\phi \in \kappa.\tau}
    \left(\theta_{\phi}^{[\beta]\kappa, \tau}\right)^{2}
    C_{\phi}^{\beta}(\mathbf{X}).
\end{equation}
For constant $a$ and $b$
\begin{equation}\label{pi8}
    C_{\phi}^{[\beta]\kappa, \tau}(a\mathbf{X}, b\mathbf{X}) = a^{k} b^{t} \
    C_{\phi}^{[\beta]\kappa, \tau}(\mathbf{X}, \mathbf{Y}).
\end{equation}
The next expansion can be used to derive several useful results of invariant
polynomials. From (\ref{pi1}), (\ref{pi5}) and (\ref{fh00}) we obtain
\begin{eqnarray}
  \int_{\mathbf{H} \in \mathfrak{U}^{\beta}(m)}  \etr\{\mathbf{AH}^{*}\mathbf{XH}
  + \mathbf{BH}^{*}\mathbf{YH}\}(d\mathbf{H}) \hspace{4cm} \nonumber\\\label{pi9}
  = \sum_{\kappa,\tau;\phi}^{\infty}
  \frac{C_{\phi}^{[\beta]\kappa, \tau}(\mathbf{A}, \mathbf{B})
  C_{\phi}^{[\beta]\kappa, \tau}(\mathbf{X},
  \mathbf{Y})}{k!t!C_{\phi}^{\beta}(\mathbf{I})},
\end{eqnarray}
where
$$
  \sum_{\kappa,\tau;\phi}^{\infty} = \sum_{k=0}^{\infty}\sum_{t=0}^{\infty}\sum_{\kappa}\sum_{\tau}\sum_{\phi \in
  \kappa.\tau}.
$$
From (\ref{pi9}) we obtain, see \citet{dg:08},
\begin{eqnarray}
    \int_{\mathbf{H} \in \mathfrak{U}^{\beta}(m)}
    C_{\phi}^{[\beta]\kappa, \tau}(\mathbf{AH}^{*}\mathbf{XH},
    \mathbf{BH}^{*}\mathbf{YH})(d\mathbf{H}) \hspace{4cm}\nonumber\\\label{pi10}
    = \frac{C_{\phi}^{[\beta]\kappa, \tau}(\mathbf{A}, \mathbf{B})
    C_{\phi}^{[\beta]\kappa, \tau}(\mathbf{X},
    \mathbf{Y})}{\theta_{\phi}^{[\beta]\kappa, \tau}C_{\phi}^{\beta}(\mathbf{I})},
\end{eqnarray}
in particular using (\ref{pi3})
\begin{equation}\label{pi11}
    \int_{\mathbf{H} \in \mathfrak{U}^{\beta}(m)}
    C_{\phi}^{[\beta]\kappa, \tau}(\mathbf{A}^{*}\mathbf{H}^{*}\mathbf{XHA},
    \mathbf{B})(d\mathbf{H}) =
    \frac{C_{\phi}^{[\beta]\kappa, \tau}(\mathbf{A}^{*}\mathbf{A}, \mathbf{B})
    C_{\kappa}^{\beta}(\mathbf{X})} {C_{\phi}^{\beta}(\mathbf{I})},
\end{equation}
analogously
\begin{equation}\label{pi12}
    \int_{\mathbf{H} \in \mathfrak{U}^{\beta}(m)}
    C_{\phi}^{[\beta]\kappa, \tau}(\mathbf{A},
    \mathbf{B}^{*}\mathbf{H}^{*}\mathbf{YHB})(d\mathbf{H}) =
    \frac{C_{\phi}^{[\beta]\kappa, \tau}(\mathbf{A}, \mathbf{B}^{*}\mathbf{B})
    C_{\tau}^{\beta}(\mathbf{Y})}{C_{\phi}^{\beta}(\mathbf{I})},
\end{equation}

\bigskip%
\textbf{Laplace transform}.%
\bigskip

For all $\mathbf{A}$ and $\mathbf{B} \in \mathfrak{S}_{m}^{\beta}$, $\mathbf{Z} \in
\mathbf{\Phi}$
\begin{eqnarray}
  \int_{\mathbf{X} \in \mathfrak{P}_{m}^{\beta}} \etr\{-\mathbf{XZ}\}|\mathbf{X}|^{a-(m-1)\beta/2-1}
    C_{\phi}^{[\beta]\kappa, \tau}(\mathbf{AX},\mathbf{BX})(d\mathbf{X}) \hspace{2.5cm} \nonumber\\ \label{pi14}
    = \Gamma_{m}^{\beta}[a, \phi]|\mathbf{Z}|^{-a}
    C_{\phi}^{[\beta]\kappa, \tau}(\mathbf{AZ}^{-1},\mathbf{BZ}^{-1}).
\end{eqnarray}
valid for $\re(a) > (m-1)\beta/2 + (k+t)_{1}$.  In particular
\begin{eqnarray}
  \int_{\mathbf{X} \in \mathfrak{P}_{m}^{\beta}} \etr\{-\mathbf{XZ}\}|\mathbf{X}|^{a-(m-1)\beta/2-1}
     C_{\phi}^{[\beta]\kappa, \tau}(\mathbf{AXA}^{*},
    \mathbf{B})(d\mathbf{X}) \hspace{2cm} \nonumber\\ \label{pi15}
    = \Gamma_{m}^{\beta}[a, \kappa]|\mathbf{Z}|^{-a}
    C_{\phi}^{[\beta]\kappa, \tau}(\mathbf{AZ}^{-1}\mathbf{A}^{*},\mathbf{B}),
\end{eqnarray}
where $\re(a) > (m-1)\beta/2 + k_{1}$. And
\begin{eqnarray}
  \int_{\mathbf{X} \in \mathfrak{P}_{m}^{\beta}} \etr\{-\mathbf{XZ}\}|\mathbf{X}|^{a-(m-1)\beta/2-1}
     C_{\phi}^{[\beta]\kappa, \tau}(\mathbf{A},
    \mathbf{BXB}^{*})(d\mathbf{X}) \hspace{2cm} \nonumber\\ \label{pi16}
    = \Gamma_{m}^{\beta}[a, \tau]|\mathbf{Z}|^{-a}
    C_{\phi}^{[\beta]\kappa, \tau}(\mathbf{A},\mathbf{BZ}^{-1}\mathbf{B}^{*}),
\end{eqnarray}
with $\re(a) > (m-1)\beta/2 + t_{1}$.

Similarly, for all $\mathbf{A}$ and $\mathbf{B} \in \mathfrak{S}_{m}^{\beta}$,
$\mathbf{Z} \in \mathbf{\Phi}$,
\begin{eqnarray}
  \int_{\mathbf{X} \in \mathfrak{P}_{m}^{\beta}} \etr\{-\mathbf{XZ}\}|\mathbf{X}|^{a-(m-1)\beta/2-1}
     C_{\phi}^{[\beta]\kappa, \tau}(\mathbf{AX}^{-1},
    \mathbf{BX}^{-1})(d\mathbf{X}) \hspace{2cm} \nonumber\\ \label{pi17}
    = \Gamma_{m}^{\beta}[a, -\phi]|\mathbf{Z}|^{-a}
    C_{\phi}^{[\beta]\kappa, \tau}(\mathbf{AZ},\mathbf{BZ}),
\end{eqnarray}
where $\re(a) > (m-1)\beta/2 - (k+t)_{m}$. In particular
\begin{eqnarray}
  \int_{\mathbf{X} \in \mathfrak{P}_{m}^{\beta}} \etr\{-\mathbf{XZ}\}|\mathbf{X}|^{a-(m-1)\beta/2-1}
     C_{\phi}^{[\beta]\kappa, \tau}(\mathbf{AX}^{-1}\mathbf{A}^{*},
    \mathbf{B})(d\mathbf{X}) \hspace{2cm} \nonumber\\ \label{pi18}
    = \Gamma_{m}^{\beta}[a, -\kappa]|\mathbf{Z}|^{-a}
    C_{\phi}^{[\beta]\kappa, \tau}(\mathbf{AZA}^{*},\mathbf{B}),
\end{eqnarray}
with $\re(a) > (m-1)\beta/2 -k_{m}$. And
\begin{eqnarray}
  \int_{\mathbf{X} \in \mathfrak{P}_{m}^{\beta}} \etr\{-\mathbf{XZ}\}|\mathbf{X}|^{a-(m-1)\beta/2-1}
     C_{\phi}^{[\beta]\kappa, \tau}(\mathbf{A},
    \mathbf{BX}^{-1}\mathbf{B}^{*})(d\mathbf{X}) \hspace{2cm} \nonumber\\ \label{pi19}
    = \Gamma_{m}^{\beta}[a, -\tau]|\mathbf{Z}|^{-a}
    C_{\phi}^{[\beta]\kappa, \tau}(\mathbf{A},\mathbf{BZB}^{*}),
\end{eqnarray}
valid  for $\re(a) > (m-1)\beta/2 -t_{m}$.

\bigskip%
\textbf{Inverse Laplace transform}.%
\bigskip

Assume that $\mathbf{Z}$ and $\mathbf{X} \in \mathfrak{S}_{m}^{\beta, \mathfrak{C}}$,
$\mathbf{A}$ and $\mathbf{B} \in \mathfrak{S}_{m}^{\beta}$ and $\re(b) > b_{0}$. Then
\begin{eqnarray}
  \frac{\Gamma_{m}^{\beta}[b, \phi] 2^{m(m-1)\beta/2}}{(2\pi i)^{m(m-1)\beta/2+ m}} \int_{\mathbf{Z}
  - \mathbf{Z}_{0} \in \mathbf{\Phi}} \etr\{\mathbf{XZ}\} |\mathbf{Z}|^{-b}
    C_{\phi}^{[\beta]\kappa, \tau}(\mathbf{AZ}^{-1},\mathbf{BZ}^{-1})(d\mathbf{Z}) \nonumber\\\label{piilt1}
    \hspace{3cm}= |\mathbf{X}|^{b-(m-1)\beta/2-1}C_{\phi}^{[\beta]\kappa, \tau}(\mathbf{AX},\mathbf{BX}),
\end{eqnarray}
and
\begin{eqnarray}
  \frac{\Gamma_{m}^{\beta}[b, -\phi] 2^{m(m-1)\beta/2}}{(2\pi i)^{m(m-1)\beta/2+ m}} \int_{\mathbf{Z}
  - \mathbf{Z}_{0} \in \mathbf{\Phi}} \etr\{\mathbf{XZ}\} |\mathbf{Z}|^{-b}
    C_{\phi}^{[\beta]\kappa, \tau}(\mathbf{AZ},\mathbf{BZ})(d\mathbf{Z}) \hspace{0.5cm}\nonumber\\\label{piilt2}
    \hspace{3cm}= |\mathbf{X}|^{b-(m-1)\beta/2-1}C_{\phi}^{[\beta]\kappa, \tau}(\mathbf{AX}^{-1},
    \mathbf{BX}^{-1}).
\end{eqnarray}
Similar expressions are obtained for
$C_{\phi}^{[\beta]\kappa,\tau}(\mathbf{AZ}^{-1},\mathbf{B})$ and
$C_{\phi}^{[\beta]\kappa, \tau}(\mathbf{A},\mathbf{BZ}^{-1})$ from (\ref{piilt1});
and for $C_{\phi}^{[\beta]\kappa, \tau}(\mathbf{AZ},\mathbf{B})$ and
$C_{\phi}^{[\beta]\kappa, \tau}(\mathbf{A},\mathbf{BZ})$ from (\ref{piilt2}).

\bigskip%
\textbf{Beta type I integrals}.%
\bigskip

For all $\mathbf{A}$ and $\mathbf{B} \in \mathfrak{S}_{m}^{\beta,\mathfrak{C}}$  and
$\re(b) > (m-1)\beta/2$,
\begin{eqnarray}
  \int_{\mathbf{0} < \mathbf{X} < \mathbf{I}} |\mathbf{X}|^{a-(m-1)\beta/2-1} |\mathbf{I} -
    \mathbf{X}|^{b-(m-1)\beta/2-1} C_{\phi}^{[\beta]\kappa, \tau}(\mathbf{AX},\mathbf{BX})(d\mathbf{X})
    \hspace{0.5cm} \nonumber\\\label{pi20}
    = \frac{\Gamma_{m}^{\beta}[a, \phi]\Gamma_{m}^{\beta}[b]}{\Gamma_{m}^{\beta}[a + b, \phi]}
    C_{\phi}^{[\beta]\kappa, \tau}(\mathbf{A},\mathbf{B}),
\end{eqnarray}
valid for $\re(a) > (m-1)\beta/2 - (k+t)_{m}$. In particular
\begin{eqnarray}
  \int_{\mathbf{0} < \mathbf{X} < \mathbf{I}} |\mathbf{X}|^{a-(m-1)\beta/2-1} |\mathbf{I} -
    \mathbf{X}|^{b-(m-1)\beta/2-1} C_{\phi}^{[\beta]\kappa, \tau}(\mathbf{AXA}^{*},\mathbf{B})(d\mathbf{X})
    \hspace{.5cm} \nonumber\\\label{pi21}
    = \frac{\Gamma_{m}^{\beta}[a, \kappa]\Gamma_{m}^{\beta}[b]}{\Gamma_{m}^{\beta}[a + b, \kappa]}
    C_{\phi}^{[\beta]\kappa, \tau}(\mathbf{AA}^{*},\mathbf{B}),
\end{eqnarray}
with $\re(a) > (m-1)\beta/2 - k_{m}$. And
\begin{eqnarray}
  \int_{\mathbf{0} < \mathbf{X} < \mathbf{I}} |\mathbf{X}|^{a-(m-1)\beta/2-1} |\mathbf{I} -
    \mathbf{X}|^{b-(m-1)\beta/2-1} C_{\phi}^{[\beta]\kappa, \tau}(\mathbf{A},\mathbf{BXB}^{*})(d\mathbf{X})
    \hspace{.5cm} \nonumber\\\label{pi22}
    = \frac{\Gamma_{m}^{\beta}[a, \tau]\Gamma_{m}^{\beta}[b]}{\Gamma_{m}^{\beta}[a + b, \tau]}
    C_{\phi}^{[\beta]\kappa, \tau}(\mathbf{A},\mathbf{BB}^{*}).
\end{eqnarray}
where $\re(a) > (m-1)\beta/2 - t_{m}$. Another particular integral given in the real
case by \citet{da:79} is
\begin{eqnarray}
  \int_{\mathbf{0} < \mathbf{X} < \mathbf{I}} |\mathbf{X}|^{a-(m-1)\beta/2-1} |\mathbf{I} -
    \mathbf{X}|^{b-(m-1)\beta/2-1} C_{\phi}^{[\beta]\kappa, \tau}(\mathbf{X}, \mathbf{I} - \mathbf{X}^{*})(d\mathbf{X})
    \hspace{.5cm} \nonumber\\\label{pi23}
    = \frac{\Gamma_{m}^{\beta}[a, \kappa]\Gamma_{m}^{\beta}[b,\tau]}{\Gamma_{m}^{\beta}[a + b, \phi]}
    \theta_{\phi}^{[\beta]\kappa, \tau} C_{\phi}^{\beta}(\mathbf{I}),
\end{eqnarray}
valid for $\re(a) > (m-1)\beta/2 - k_{m}$ and $\re(b) > (m-1)\beta/2 - k_{m}$.

Analogously,  for all $\mathbf{A}$ and $\mathbf{B} \in
\mathfrak{S}_{m}^{\beta,\mathfrak{C}}$ and $\re(b) > (m-1)\beta/2$,
\begin{eqnarray}
  \int_{\mathbf{0} < \mathbf{X} < \mathbf{I}} |\mathbf{X}|^{a-(m-1)\beta/2-1} |\mathbf{I} -
    \mathbf{X}|^{b-(m-1)\beta/2-1} C_{\phi}^{[\beta]\kappa, \tau}(\mathbf{AX}^{-1},\mathbf{BX}^{-1})(d\mathbf{X})
    \nonumber\\\label{pi24}
    = \frac{\Gamma_{m}^{\beta}[a, -\phi]\Gamma_{m}^{\beta}[b]}{\Gamma_{m}^{\beta}[a + b, -\phi]}
    C_{\phi}^{[\beta]\kappa, \tau}(\mathbf{A},\mathbf{B}),
\end{eqnarray}
where  $\re(a) > (m-1)\beta/2 + (k+t)_{1}$. In particular
\begin{eqnarray}
  \int_{\mathbf{0} < \mathbf{X} < \mathbf{I}} |\mathbf{X}|^{a-(m-1)\beta/2-1} |\mathbf{I} -
    \mathbf{X}|^{b-(m-1)\beta/2-1} C_{\phi}^{[\beta]\kappa, \tau}(\mathbf{AX}^{-1}\mathbf{A}^{*},\mathbf{B})(d\mathbf{X})
    \nonumber\\\label{pi25}
    = \frac{\Gamma_{m}^{\beta}[a, -\kappa]\Gamma_{m}^{\beta}[b]}{\Gamma_{m}^{\beta}[a + b, -\kappa]}
    C_{\phi}^{[\beta]\kappa, \tau}(\mathbf{AA}^{*},\mathbf{B}),
\end{eqnarray}
valid for  $\re(a) > (m-1)\beta/2 + k_{1}$. And
\begin{eqnarray}
  \int_{\mathbf{0} < \mathbf{X} < \mathbf{I}} |\mathbf{X}|^{a-(m-1)\beta/2-1} |\mathbf{I} -
    \mathbf{X}|^{b-(m-1)\beta/2-1} C_{\phi}^{[\beta]\kappa, \tau}(\mathbf{A},\mathbf{BX}^{-1}\mathbf{B}^{*})(d\mathbf{X})
    \nonumber\\\label{pi26}
    = \frac{\Gamma_{m}^{\beta}[a, -\tau]\Gamma_{m}^{\beta}[b]}{\Gamma_{m}^{\beta}[a + b, -\tau]}
    C_{\phi}^{[\beta]\kappa, \tau}(\mathbf{A},\mathbf{BB}^{*}),
\end{eqnarray}
with $\re(a) > (m-1)\beta/2 + t_{1}$.

Now, taking $b = (m-1)\beta/2+1 > (m-1)\beta/2$ in (\ref{pi20}) and (\ref{pi24}) we
have the following results.

For all $\mathbf{A}$ and $\mathbf{B} \in \mathfrak{S}_{m}^{\beta,\mathfrak{C}}$,
\begin{eqnarray}
  \int_{\mathbf{0} < \mathbf{X} < \mathbf{I}} |\mathbf{X}|^{a-(m-1)\beta/2-1}
    C_{\phi}^{[\beta]\kappa, \tau}(\mathbf{AX},\mathbf{BX})(d\mathbf{X})
    \hspace{3.5cm} \nonumber\\\label{pi26a}
    = \frac{\Gamma_{m}^{\beta}[a, \phi]\Gamma_{m}^{\beta}[(m-1)\beta/2+1]}
    {\Gamma_{m}^{\beta}[a + (m-1)\beta/2+1, \phi]}
    C_{\phi}^{[\beta]\kappa, \tau}(\mathbf{A},\mathbf{B}),
\end{eqnarray}
valid for $\re(a) > (m-1)\beta/2 - (k+t)_{m}$. And,
\begin{eqnarray}
  \int_{\mathbf{0} < \mathbf{X} < \mathbf{I}} |\mathbf{X}|^{a-(m-1)\beta/2-1}
    C_{\phi}^{[\beta]\kappa, \tau}(\mathbf{AX}^{-1},\mathbf{BX}^{-1})(d\mathbf{X})
    \hspace{3cm} \nonumber\\\label{pi26b}
    = \frac{\Gamma_{m}^{\beta}[a, -\phi]\Gamma_{m}^{\beta}[(m-1)\beta/2+1]}
    {\Gamma_{m}^{\beta}[a + (m-1)\beta/2+1, -\phi]}
    C_{\phi}^{[\beta]\kappa, \tau}(\mathbf{A},\mathbf{B}),
\end{eqnarray}
where  $\re(a) > (m-1)\beta/2 + (k+t)_{1}$.

\bigskip%
\textbf{Beta type II integrals}.%
\bigskip

For all $\mathbf{A}$ and $\mathbf{B} \in \mathfrak{S}_{m}^{\beta,\mathfrak{C}}$,
\begin{eqnarray}
  \int_{\mathbf{X} \in \mathfrak{P}_{m}^{\beta}} |\mathbf{X}|^{a-(m-1)\beta/2-1}
  |\mathbf{I} + \mathbf{X}|^{-(a+b)} C_{\phi}^{[\beta]\kappa, \tau}(\mathbf{AX},\mathbf{BX})
  \hspace{2cm} \nonumber\\\label{pi27}
  = \frac{\Gamma_{m}^{\beta}[a, \phi] \Gamma_{m}^{\beta}[b, -\phi]}{\Gamma_{m}^{\beta}[a + b]}
  C_{\phi}^{[\beta]\kappa, \tau}(\mathbf{A},\mathbf{B}),
\end{eqnarray}
with $\re(a) > (m-1)\beta/2 - (k+t)_{m}$ and $\re(b) > (m-1)\beta/2 + (k+t)_{1}$. In
particular
\begin{eqnarray}
  \int_{\mathbf{X} \in \mathfrak{P}_{m}^{\beta}} |\mathbf{X}|^{a-(m-1)\beta/2-1}
  |\mathbf{I} + \mathbf{X}|^{-(a+b)} C_{\phi}^{[\beta]\kappa, \tau}(\mathbf{AXA}^{*},\mathbf{B})
  \hspace{2cm} \nonumber\\\label{pi28}
  = \frac{\Gamma_{m}^{\beta}[a, \kappa] \Gamma_{m}^{\beta}[b, -\kappa]}{\Gamma_{m}^{\beta}[a + b]}
  C_{\phi}^{[\beta]\kappa, \tau}(\mathbf{AA}^{*},\mathbf{B}),
\end{eqnarray}
such that $\re(a) > (m-1)\beta/2 - k_{m}$ and $\re(b) > (m-1)\beta/2 + k_{1}$. And
\begin{eqnarray}
  \int_{\mathbf{X} \in \mathfrak{P}_{m}^{\beta}} |\mathbf{X}|^{a-(m-1)\beta/2-1}
  |\mathbf{I} + \mathbf{X}|^{-(a+b)} C_{\phi}^{[\beta]\kappa, \tau}(\mathbf{A},\mathbf{BXB}^{*})
  \hspace{2cm} \nonumber\\\label{pi29}
  = \frac{\Gamma_{m}^{\beta}[a, \tau] \Gamma_{m}^{\beta}[b, -\tau]}{\Gamma_{m}^{\beta}[a + b]}
  C_{\phi}^{[\beta]\kappa, \tau}(\mathbf{A},\mathbf{BB}^{*}),
\end{eqnarray}
valid for $\re(a) > (m-1)\beta/2 - t_{m}$ and $\re(b) > (m-1)\beta/2 + t_{1}$.

In a similar way, for all $\mathbf{A}$ and $\mathbf{B} \in
\mathfrak{S}_{m}^{\beta,\mathfrak{C}}$,
\begin{eqnarray}
  \int_{\mathbf{X} \in \mathfrak{P}_{m}^{\beta}} |\mathbf{X}|^{a-(m-1)\beta/2-1}
  |\mathbf{I} + \mathbf{X}|^{-(a+b)} C_{\phi}^{[\beta]\kappa, \tau}(\mathbf{AX}^{-1},\mathbf{BX}^{-1})
  \hspace{2cm} \nonumber\\\label{pi30}
  = \frac{\Gamma_{m}^{\beta}[a, -\phi] \Gamma_{m}^{\beta}[b, \phi]}{\Gamma_{m}^{\beta}[a + b]}
  C_{\phi}^{[\beta]\kappa, \tau}(\mathbf{A},\mathbf{B})
\end{eqnarray}
where $\re(a) > (m-1)\beta/2 + (k+t)_{1}$ and $\re(b) > (m-1)\beta/2 - (k+t)_{m}$. In
particular
\begin{eqnarray}
  \int_{\mathbf{X} \in \mathfrak{P}_{m}^{\beta}} |\mathbf{X}|^{a-(m-1)\beta/2-1}
  |\mathbf{I} + \mathbf{X}|^{-(a+b)} C_{\phi}^{[\beta]\kappa, \tau}(\mathbf{AX}^{-1}\mathbf{A}^{*},\mathbf{B})
  \hspace{2cm} \nonumber\\\label{pi31}
  = \frac{\Gamma_{m}^{\beta}[a, -\kappa] \Gamma_{m}^{\beta}[b, \kappa]}{\Gamma_{m}^{\beta}[a + b]}
  C_{\phi}^{[\beta]\kappa, \tau}(\mathbf{AA}^{*},\mathbf{B}),
\end{eqnarray}
with $\re(a) > (m-1)\beta/2 + k_{1}$ and $\re(b) > (m-1)\beta/2 - k_{m}$. And
\begin{eqnarray}
  \int_{\mathbf{X} \in \mathfrak{P}_{m}^{\beta}} |\mathbf{X}|^{a-(m-1)\beta/2-1}
  |\mathbf{I} + \mathbf{X}|^{-(a+b)} C_{\phi}^{[\beta]\kappa, \tau}(\mathbf{A},\mathbf{BX}^{-1}\mathbf{B}^{*})
  \hspace{2cm} \nonumber\\\label{pi32}
  = \frac{\Gamma_{m}^{\beta}[a, -\tau] \Gamma_{m}^{\beta}[b, \tau]}{\Gamma_{m}^{\beta}[a + b]}
  C_{\phi}^{[\beta]\kappa, \tau}(\mathbf{A},\mathbf{BB}^{*}),
\end{eqnarray}
such that $\re(a) > (m-1)\beta/2 + t_{1}$ and $\re(b) > (m-1)\beta/2 - t_{m}$.

\bigskip%
\textbf{Incomplete gamma and beta functions}.%
\bigskip

First consider the following results

For all $\mathbf{A}$ and $\mathbf{B} \in \mathfrak{S}_{m}^{\beta}$ and $\mathbf{0} <
\mathbf{\Xi} < \mathbf{I}$,
\begin{eqnarray}
  \int_{\mathbf{0} < \mathbf{X}< \mathbf{\Xi}} |\mathbf{X}|^{a-(m-1)\beta/2-1}
  C_{\phi}^{[\beta]\kappa, \tau}(\mathbf{AX}, \mathbf{BX})(d\mathbf{X})
  \hspace{3.5cm} \nonumber\\ \label{pi33}
  = \frac{\Gamma_{m}^{\beta}[a,\phi] \Gamma_{m}^{\beta}[(m-1)\beta/2 + 1]}
  {\Gamma_{m}^{\beta}[a+(m-1)\beta/2 + 1, \phi]}
  |\mathbf{\Omega}|^{a}C_{\phi}^{[\beta]\kappa, \tau} (\mathbf{A\Xi},\mathbf{B\Xi}).
\end{eqnarray}
valid for $\re(a) > (m-1)\beta/2 - (k+t)_{m}$. And
\begin{eqnarray}
  \int_{\mathbf{0} < \mathbf{X}< \mathbf{\Xi}} |\mathbf{X}|^{a-(m-1)\beta/2-1}
  C_{\phi}^{[\beta]\kappa, \tau}(\mathbf{AXA}^{*}, \mathbf{B})(d\mathbf{X})
  \hspace{3.5cm} \nonumber\\ \label{pi34}
  = \frac{\Gamma_{m}^{\beta}[a,\kappa] \Gamma_{m}^{\beta}[(m-1)\beta/2 + 1]}
  {\Gamma_{m}^{\beta}[a+(m-1)\beta/2 + 1, \kappa]}
  |\mathbf{\Omega}|^{a}C_{\phi}^{[\beta]\kappa, \tau} (\mathbf{A\Xi A}^{*},\mathbf{B}).
\end{eqnarray}
valid for $\re(a) > (m-1)\beta/2 - k_{m}$.

The next result is obtained immediately, expanding $\etr\{-\mathbf{XA}\}$ in terms of
Jack polynomials, making use of the property (\ref{pi5}) and integrating term by term
using (\ref{pi33}). Thus, for all $\mathbf{A}$ and $\mathbf{B} \in
\mathfrak{S}_{m}^{\beta}$ and $\mathbf{\Omega} \in \mathbf{\Phi}$,
$$
  \int_{\mathbf{0} < \mathbf{X}< \mathbf{\Omega}} \etr\{-\mathbf{XA}\}|\mathbf{X}|^{a-(m-1)\beta/2-1}
     C_{\tau}^{\beta}(\mathbf{BX})(d\mathbf{X})\hspace{5cm}
$$\vspace{-0.5cm}%
{\small
\begin{equation}\label{pi35}
    = \frac{\Gamma_{m}^{\beta}[a] \Gamma_{m}^{\beta}[(m-1)\beta/2 + 1]}{\Gamma_{m}^{\beta}[a+(m-1)\beta/2 + 1]}
    |\mathbf{\Omega}|^{a}\sum_{k = 0}^{\infty}\sum_{\kappa;\phi \in \kappa.\tau}
    \frac{[a]_{\phi}^{\beta} \theta_{\phi}^{[\beta]\kappa,\tau}C_{\phi}^{[\beta]\kappa, \tau}
    (-\mathbf{A\Omega},\mathbf{B\Omega})}{k! [a+(m-1)\beta/2 + 1]_{\phi}^{\beta}}.
\end{equation}}
valid for $\re(a) > (m-1)\beta/2 + (k+t)_{1}$.

Similarly, we expand $|\mathbf{I} - \mathbf{X}|^{b-(m-1)\beta/2-1}$ in terms of Jack
polynomials,make use of the property (\ref{pi5}) and integrate term by term using
(\ref{pi33}). For all $\mathbf{A} \in \mathfrak{S}_{m}^{\beta,\mathfrak{C}}$ and
$\mathbf{0} < \mathbf{\Xi} < \mathbf{I}$,%
{\small
$$
  \int_{\mathbf{0} < \mathbf{X}< \mathbf{\Xi}} |\mathbf{X}|^{a-(m-1)\beta/2-1} |\mathbf{I} -
    \mathbf{X}|^{b-(m-1)\beta/2-1} C_{\tau}^{\beta}(\mathbf{AX})(d\mathbf{X})
  = \frac{\Gamma_{m}^{\beta}[a] \Gamma_{m}^{\beta}[(m-1)\beta/2 + 1]}{\Gamma_{m}^{\beta}[a+(m-1)\beta/2 + 1]}
$$}
\begin{equation}\label{pi36}
    \hspace{0.5cm}\times \
    |\mathbf{\Omega}|^{a}\sum_{k = 0}^{\infty}\sum_{\kappa;\phi \in \kappa.\tau}
    \frac{[-b+(m-1)\beta/2 + 1]_{\phi}^{\beta} \theta_{\phi}^{[\beta]\kappa,\tau}
    C_{\phi}^{[\beta]\kappa, \tau}
    (\mathbf{\Omega},\mathbf{A\Omega})}{k! [a+(m-1)\beta/2 + 1]_{\phi}^{\beta}}.
\end{equation}
valid for $\re(a) > (m-1)\beta/2 - (k+t)_{m}$.

\section{Application}\label{sec6}

As an application, in this section we found the joint density eigenvalue of the
central Wishart distribution for normed division algebras, and also derived the
largest and smallest eigenvalue distributions. First, from \citet{dggj:09b} let us
consider the following definitions.

\begin{definition}
Let $\mathbf{X} \in \mathcal{L}_{m,n}^{\beta}$ be a random matrix. Then $\mathbf{X}$
is said to have a matrix variate normal distribution $\mathbf{X} \sim \mathcal{N}_{n
\times m}^{\beta}(\boldsymbol{\mu}, \mathbf{\Sigma}, \mathbf{\Theta})$, of mean
$\boldsymbol{\mu}$ and $\cov(\Tec \mathbf{X}^{T}) = \mathbf{\Theta} \otimes
\mathbf{\Sigma}$, if its density function is given by
$$
  \frac{1}{\left(2\pi \beta^{-1}\right)^{\beta mn/2}|\mathbf{\Sigma}|^{\beta n/2}
  |\mathbf{\Theta}|^{\beta m/2}}\etr\left \{- \frac{\beta}{2} \mathbf{\Sigma}^{-1}
  (\mathbf{X} - \boldsymbol{\mu})^{*} \mathbf{\Theta}^{-1}(\mathbf{X} - \boldsymbol{\mu})
  \right \}.
$$
\end{definition}
Also
\begin{definition}
Let $\mathbf{X} \in \mathcal{L}_{m,n}^{\beta}$ with distribution $\mathbf{X} \sim
\mathcal{N}_{n \times m}^{\beta}(\boldsymbol{0}, \mathbf{\Sigma}, \mathbf{I}_{n})$
and define $\mathbf{S} = \mathbf{X}^{*}\mathbf{X}$, then $\mathbf{S}$ is said to have
a central Wishart distribution $\mathbf{S} \sim \mathcal{W}_{m}^{\beta}(n,
\mathbf{\Sigma})$ with $n$  degrees of freedom and parameter $\mathbf{\Sigma}$.
Moreover, its density function is given by
$$
  \frac{1}{\left(2\beta^{-1}\right)^{\beta mn/2} \Gamma_{m}^{\beta}[\beta n/2]
  |\mathbf{\Sigma}|^{\beta n/2}} |\mathbf{S}|^{\beta(n - m +1)/2 -1}\etr\{-\beta
  \mathbf{\Sigma}^{-1}\mathbf{S}/2\}
$$
with $n \geq (m-1)\beta$.
\end{definition}
Therefore, from (\ref{dfeig}) and (\ref{exp}) the joint density of the eigenvalues,
$\lambda_{1}> \cdots > \lambda_{m} > 0$, of $\mathbf{S}$ is
\begin{eqnarray*}
  \frac{\pi^{m^{2}\beta/2 + \varrho}}{\left(2\beta^{-1}\right)^{\beta mn/2} \Gamma_{m}^{\beta}[\beta n/2]
  \Gamma_{m}^{\beta}[\beta m/2]
  |\mathbf{\Sigma}|^{\beta n/2}} \prod_{i=1}^{m}\lambda_{i}^{\beta(n - m +1)/2 -1}\hspace{3cm}\\
  \times \ \prod_{i<j}^{m}(\lambda_{i} - \lambda_{j})^{\beta} {}_{0}F_{0}(-\beta
  \mathbf{\Sigma}^{-1}/2,\mathbf{L})
\end{eqnarray*}
where $\mathbf{L} = \diag(\lambda_{1}, \dots, \lambda_{m})$.

In addition, as an immediate consequence of Theorems \ref{teoigb1} and \ref{teoigb2}
we obtain the following result.

\begin{theorem}
Let $\mathbf{S} \sim \mathcal{W}_{m}^{\beta}(n, \mathbf{\Sigma})$  and
$\mathbf{\Omega} \in \mathbf{\Phi}$, then
\begin{eqnarray}
  P(\mathbf{S} < \mathbf{\Omega})=  \frac{\Gamma_{m}^{\beta}[(m-1)\beta/2+1]}
    {\left(2\beta^{-1}\right)^{\beta mn/2} \Gamma_{m}^{\beta}[(n + m -1)\beta/2 +1]}
    \frac{|\mathbf{\Omega}|^{\beta n/2}}{|\mathbf{\Sigma}|^{\beta n/2}}\hspace{2cm}\nonumber\\ \label{fd1}
    \times \ {}_{1}F_{1}^{\beta}(\beta n/2; (n + m -1)\beta/2 +1; -\beta\mathbf{\Omega\Sigma}^{-1}/2),
\end{eqnarray}
valid for $\re(n) > (m-1)\beta - 2k_{m}$. And if $r = (n-m+1)\beta/2-1$ is a positive
integer, then
\begin{equation}\label{fd2}
    P(\mathbf{\Omega} > \mathbf{S})=   \etr\{-\beta\mathbf{\Omega\Sigma}^{-1}/2\}
   \sum_{k =0}^{mr}\sum_{\kappa}{}^{*}   \frac{C_{\kappa}^{\beta}(\beta\mathbf{\Omega\Sigma}^{-1}/2)}{k!},
\end{equation}
for $\re(n) > (m-1)\beta + 2k_{1}$ and $\sum_{\kappa}^{*}$ denotes summation over
those partitions $\kappa = (k_{1}, \dots, k_{m})$ of $k$ with $k_{1} \leq r$.
\end{theorem}

Observing that if $\lambda_{\max}$ and $\lambda_{\min}$ are the largest and smallest
eigenvalues of $\mathbf{S}$, respectively, then the inequalities $\lambda_{\max} < x$
and $\lambda_{\min} > y$ are equivalent to $\mathbf{S} < x\mathbf{I}$ and $\mathbf{S}
> y\mathbf{I}$, respectively and the following result is obtained.

\begin{corollary}
Assume that $\mathbf{S} \sim \mathcal{W}_{m}^{\beta}(n, \mathbf{\Sigma})$ and $x >
0$. Then
\begin{eqnarray}
    P(\lambda_{\max} < x)=  \frac{\Gamma_{m}^{\beta}[(m-1)\beta/2+1]}
    {\left(2\beta^{-1}\right)^{\beta mn/2} \Gamma_{m}^{\beta}[(n + m -1)\beta/2 +1]}
    \frac{x^{\beta mn/2}}{|\mathbf{\Sigma}|^{\beta n/2}}\hspace{1cm}\nonumber\\ \label{maxeig}
    \times \ {}_{1}F_{1}^{\beta}(\beta n/2; (n + m -1)\beta/2 +1; -\beta x\mathbf{\Sigma}^{-1}/2),
\end{eqnarray}
valid for $\re(n) > (m-1)\beta - 2k_{m}$. And if $r = (n-m+1)\beta/2-1$ is a positive
integer and $y > 0$, then
\begin{equation}\label{mineig}
    P(\lambda_{\min} < y)=  1 - \etr\{-\beta y \mathbf{\Sigma}^{-1}/2\}
   \sum_{k =0}^{mr}\sum_{\kappa}{}^{*}   \frac{C_{\kappa}^{\beta}(\beta y\mathbf{\Sigma}^{-1}/2)}{k!},
\end{equation}
for $\re(n) > (m-1)\beta + 2k_{1}$ and $\sum_{\kappa}^{*}$ denotes summation over
those partitions $\kappa = (k_{1}, \dots, k_{m})$ of $k$ with $k_{1} \leq r$.
\end{corollary}

As a numerical example we plot the distribution function of $\lambda_{max}$ on Figure
1 and the distribution function of $\lambda_{min}$ on Figure 2. First note that
applying the generalised Kummer relation (\ref{kr}) in (\ref{maxeig}) we obtain
\begin{eqnarray*}
    P(\lambda_{\max} < x)=  \frac{\Gamma_{m}^{\beta}[(m-1)\beta/2+1]\etr\{-\beta x \mathbf{\Sigma}^{-1}/2\}}
    {\left(2\beta^{-1}\right)^{\beta mn/2} \Gamma_{m}^{\beta}[(n + m -1)\beta/2 +1]}
    \frac{x^{\beta mn/2}}{|\mathbf{\Sigma}|^{\beta n/2}}\hspace{2cm}\nonumber\\ \label{maxeigmod}
    \times \ {}_{1}F_{1}^{\beta}((m -1)\beta/2 +1; (n + m -1)\beta/2 +1; \beta
    x\mathbf{\Sigma}^{-1}/2).
\end{eqnarray*}

\newpage

\begin{figure}[!ht]
\begin{center}
\includegraphics[angle=0,width=6cm,height=5cm]{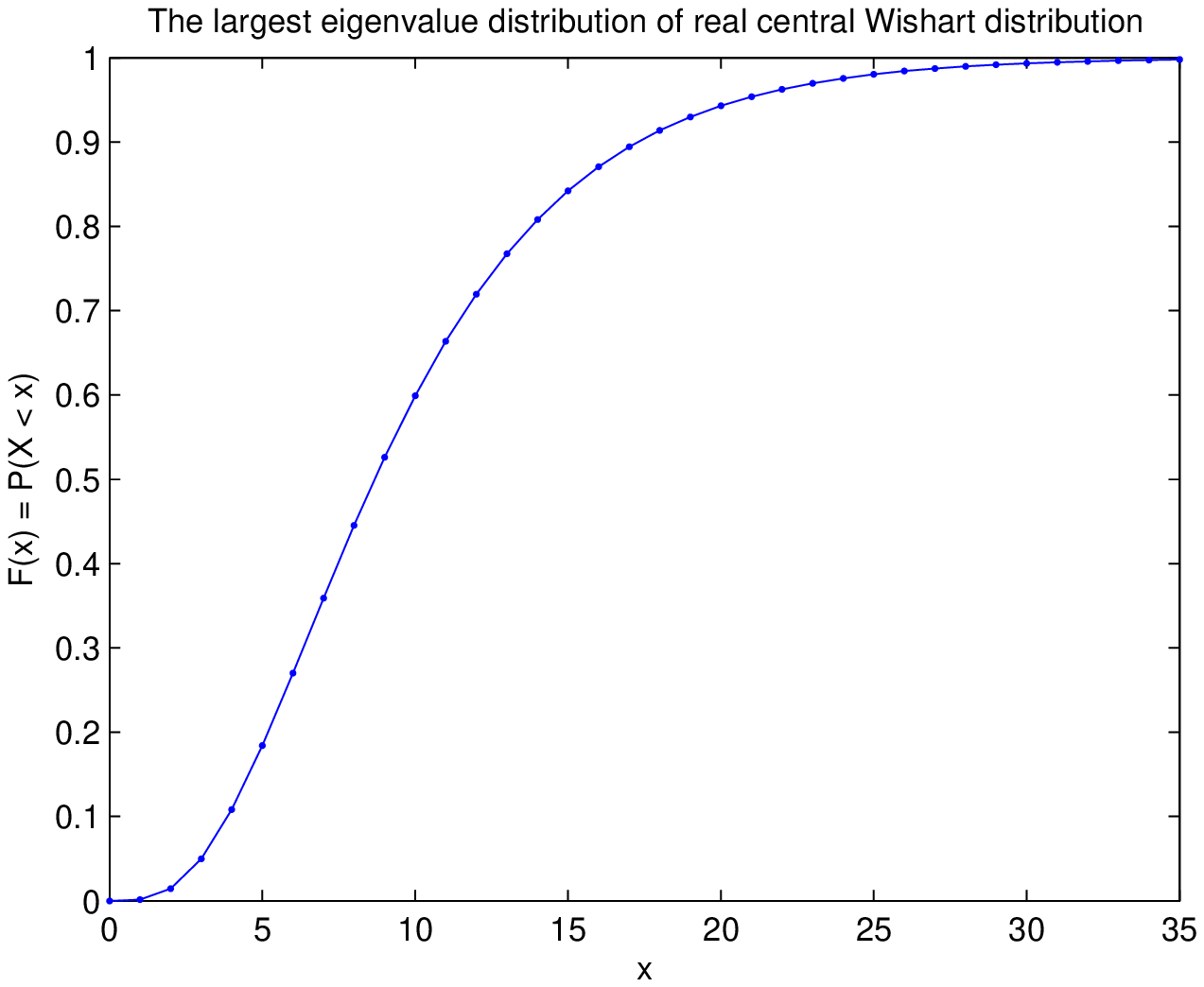}
\includegraphics[angle=0,width=6cm,height=5cm]{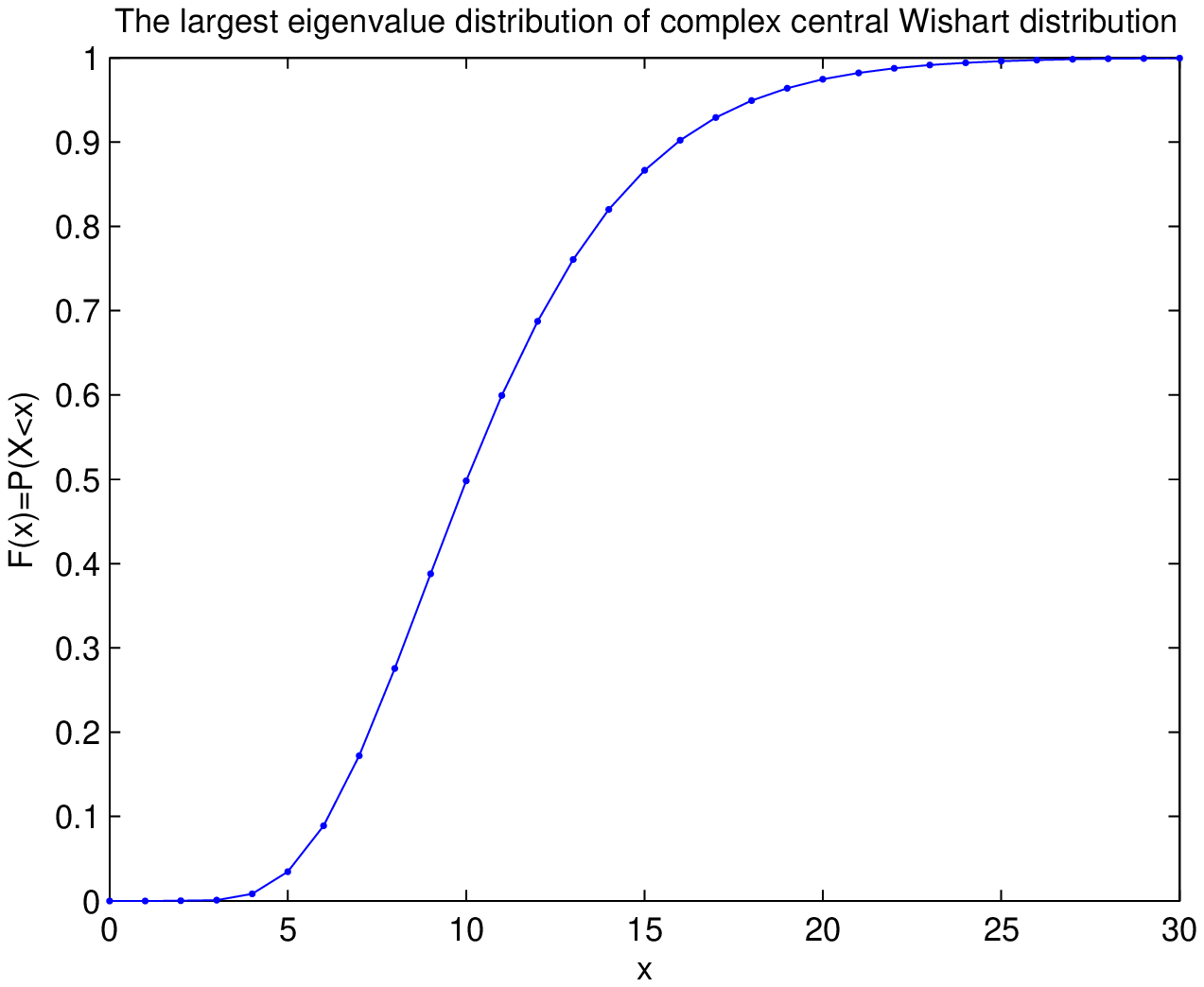}
\includegraphics[angle=0,width=6cm,height=5cm]{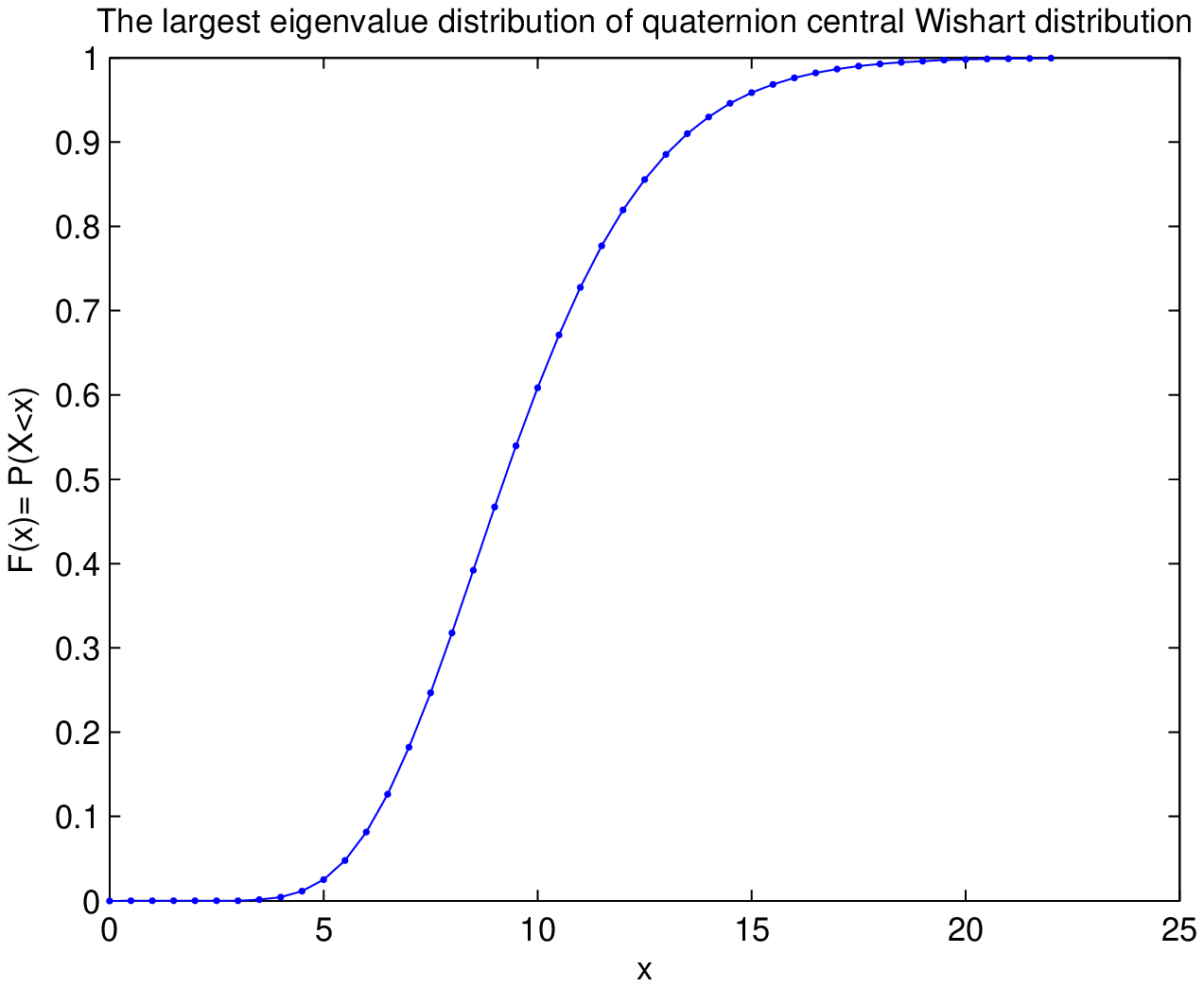}
\includegraphics[angle=0,width=6cm,height=5cm]{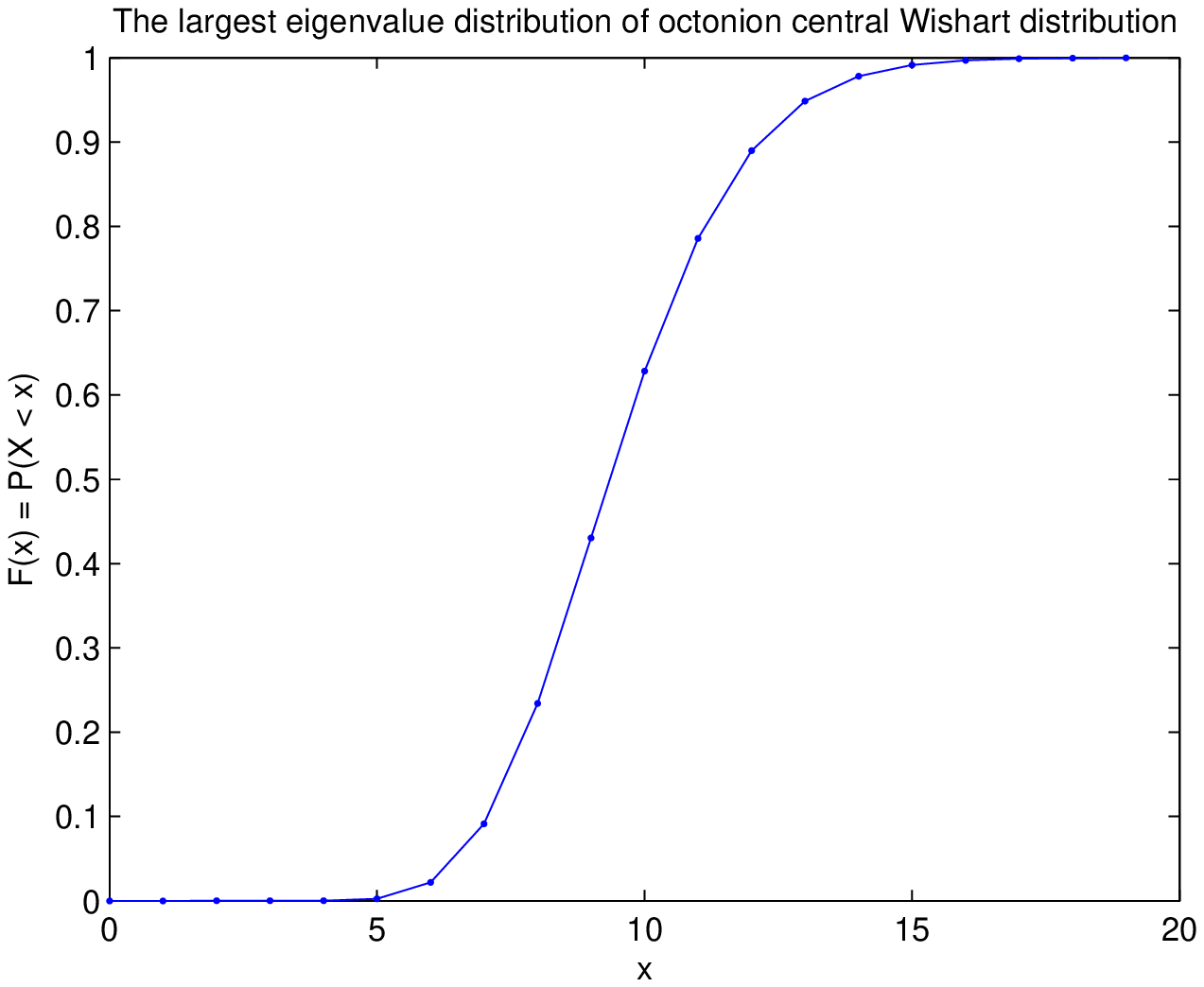}
\end{center}
\vspace{-0.9cm} \caption{Distribution functions of $\lambda_{\max}$ of
$\mathcal{W}_{2}^{\beta}(4,\diag(1,2))$, $\beta = 1, 2, 4$ and $8$.}
\end{figure}

\begin{figure}[!ht]
\begin{center}
\includegraphics[angle=0,width=6cm,height=5cm]{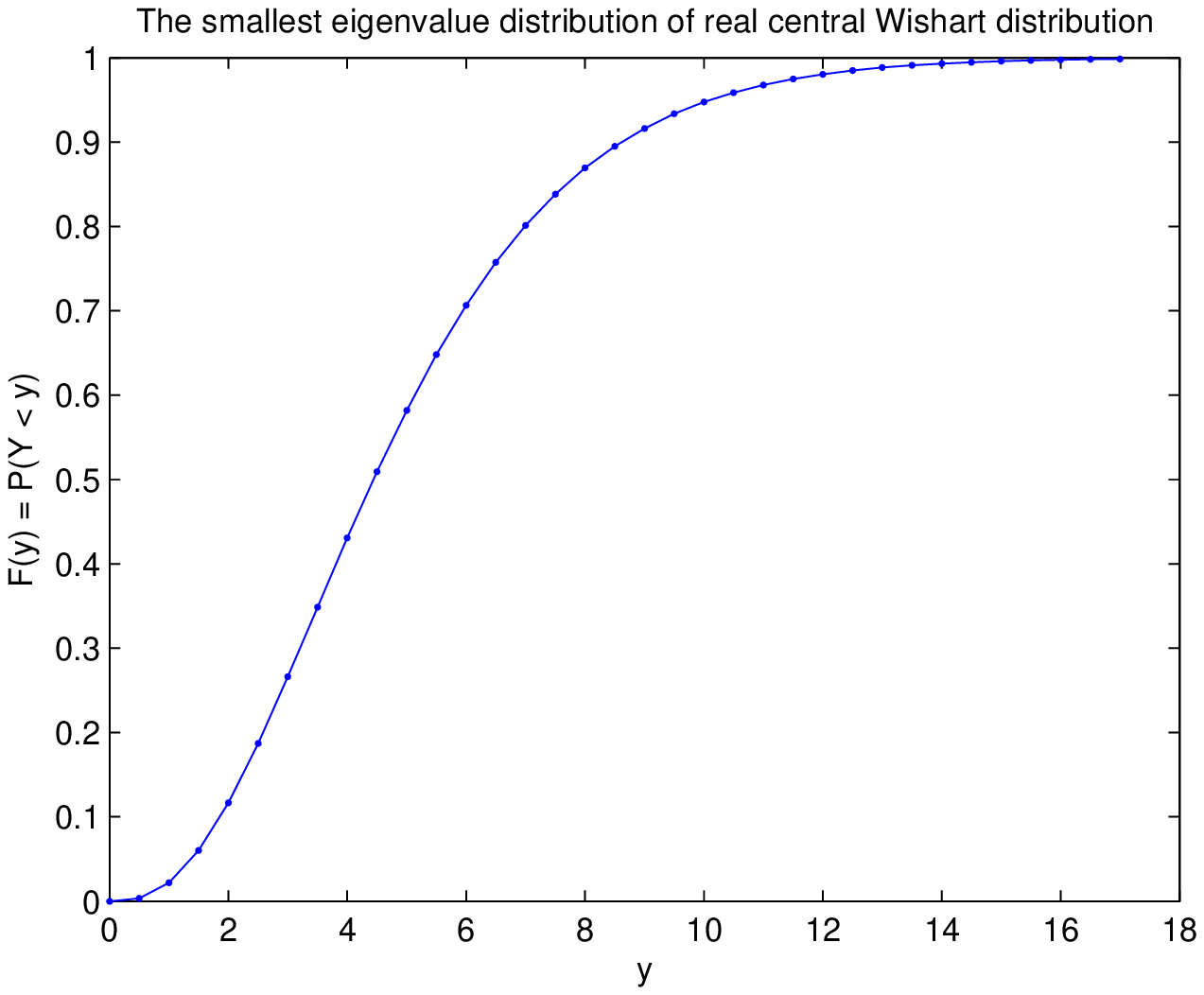}
\includegraphics[angle=0,width=6cm,height=5cm]{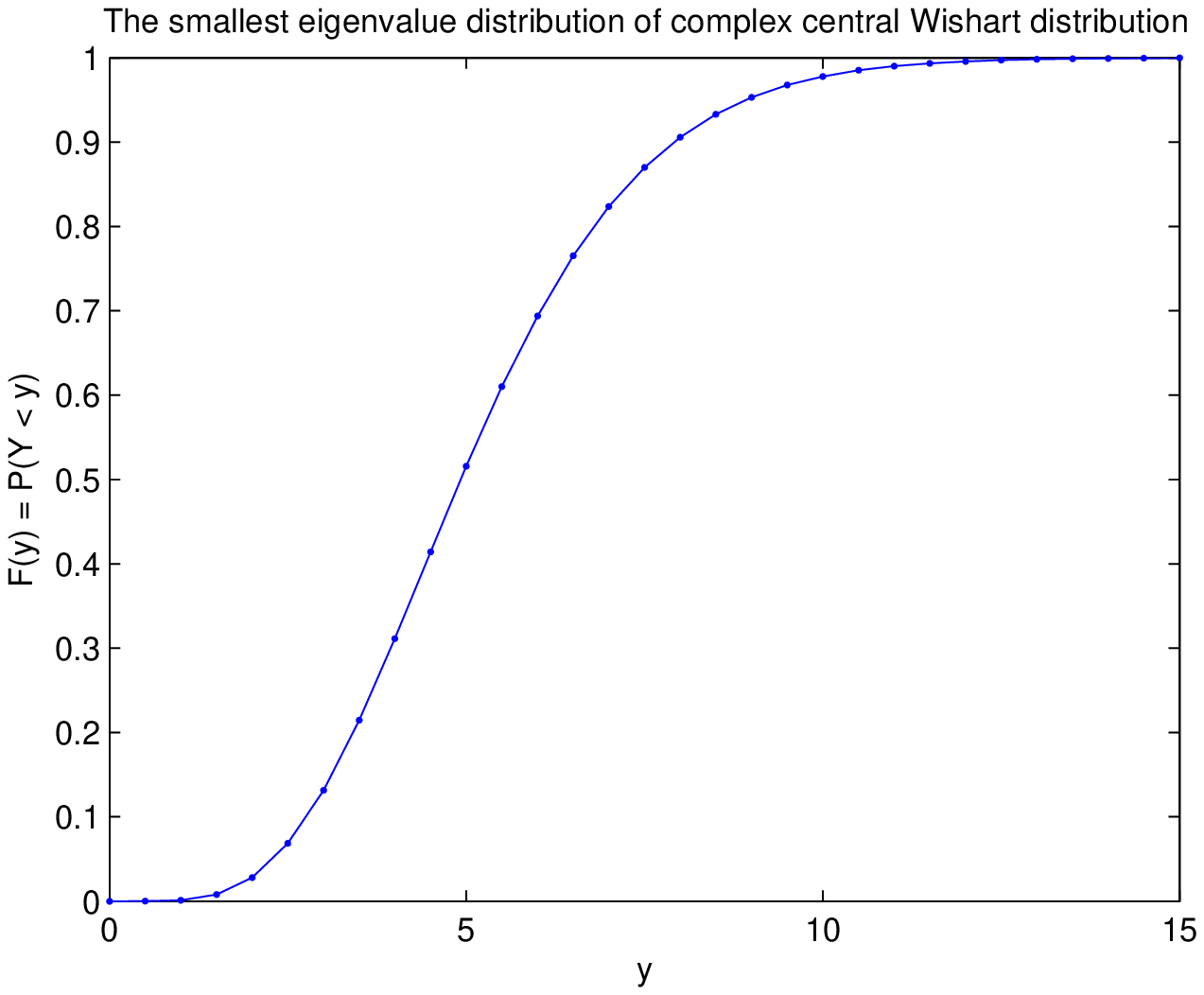}
\includegraphics[angle=0,width=6cm,height=5cm]{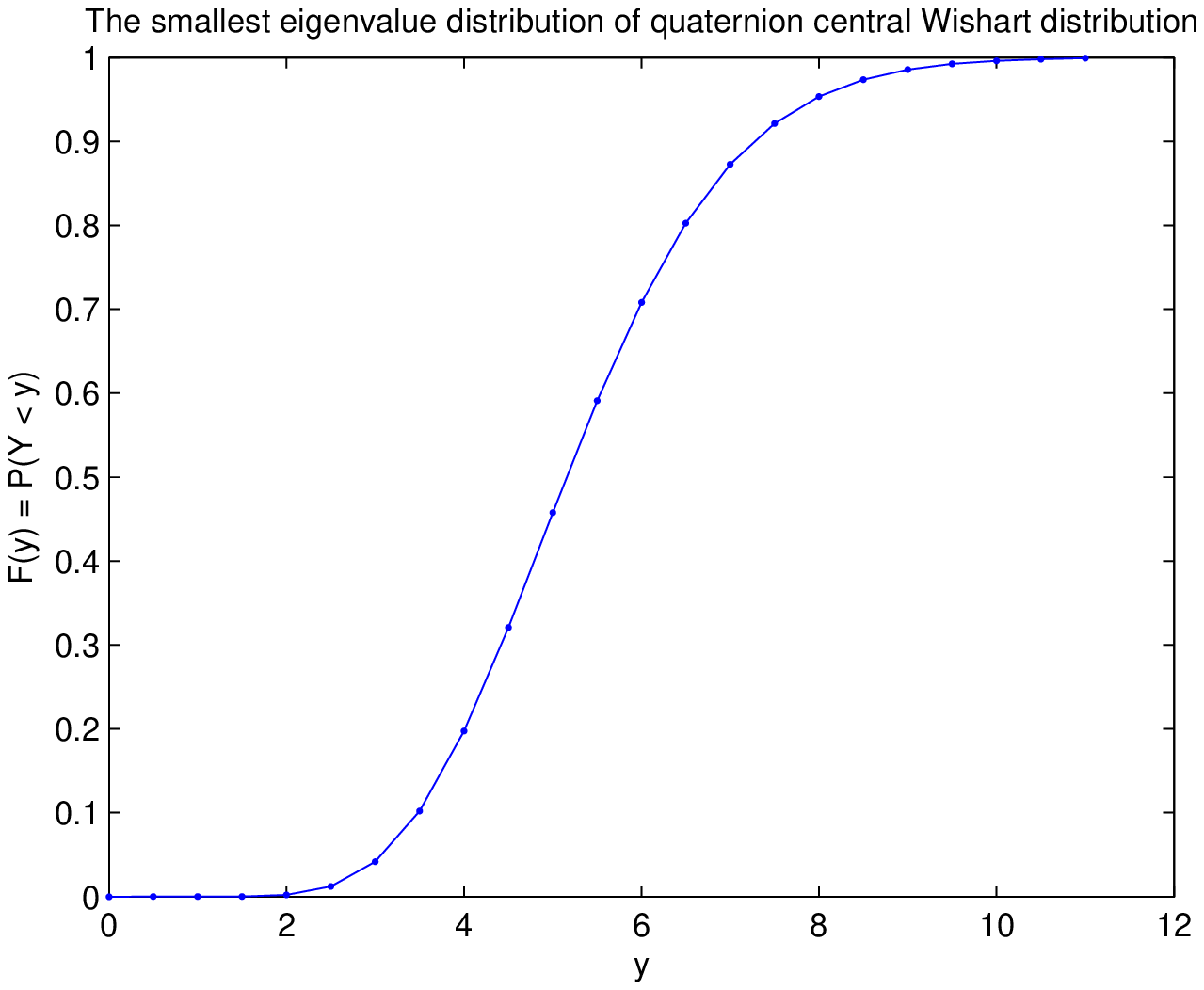}
\includegraphics[angle=0,width=6cm,height=5cm]{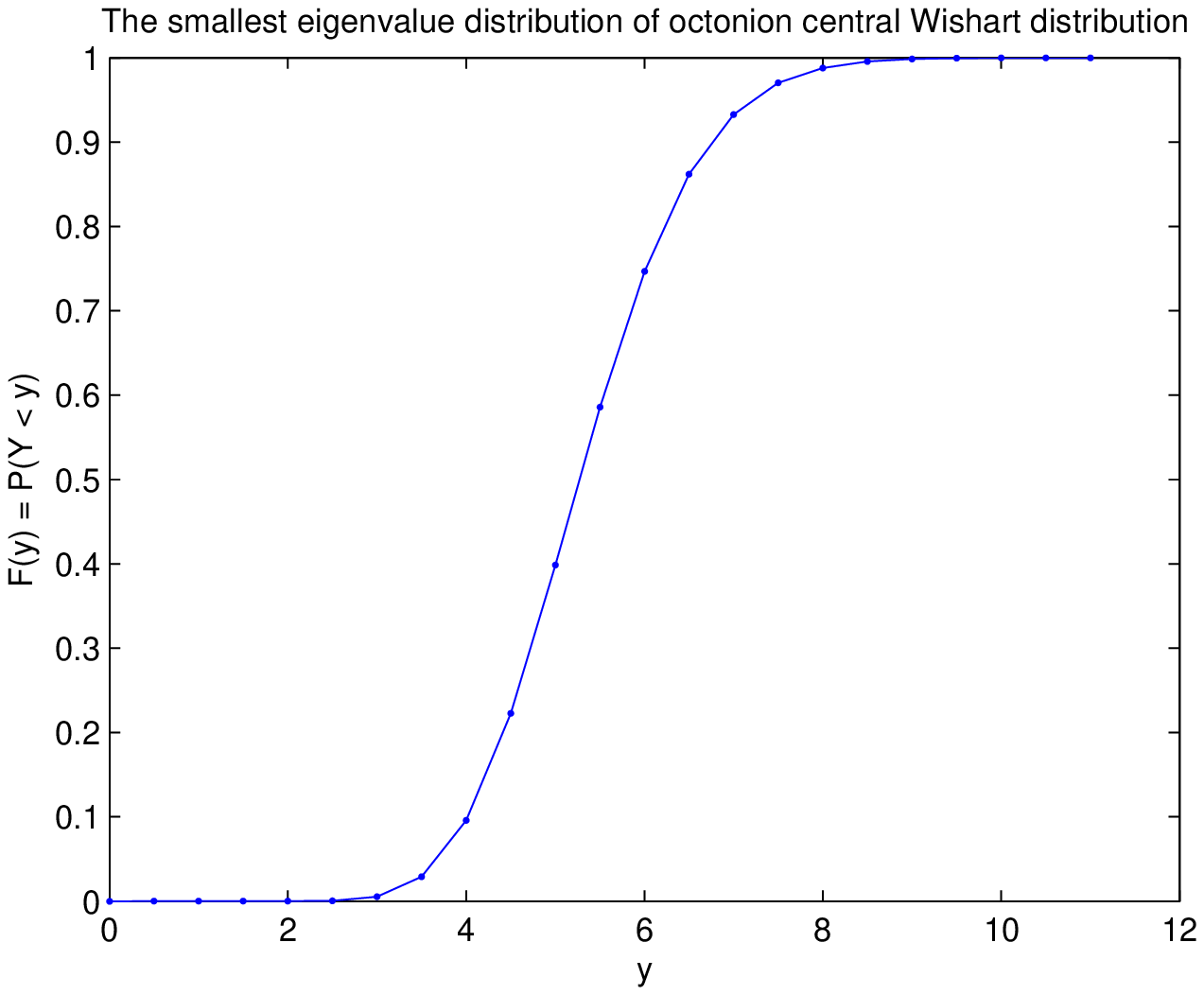}
\end{center}
\vspace{-0.9cm} \caption{Distribution functions of $\lambda_{\min}$ of
$\mathcal{W}_{2}^{\beta}(7,\diag(1,2))$, $\beta = 1, 2, 4$ and $8$.}
\end{figure}


\end{document}